\documentclass[12pt]{article}
\usepackage{amsthm,amssymb,latexsym,amsxtra,graphicx}

\theoremstyle{plain}
\newtheorem{thm}{Theorem}[section]
\newtheorem{cor}[thm]{Corollary}
\newtheorem{lem}[thm]{Lemma}
\newtheorem{mainlem}[thm]{Main Lemma}
\newtheorem{pro}[thm]{Proposition}
\newtheorem*{claim}{Claim}

\theoremstyle{definition}
\newtheorem{dfn}{Definition}
\newtheorem{block}[thm]{}

\theoremstyle{remark}
\newtheorem{rem}{Remark}
\newtheorem*{exa}{Example}
\newtheorem*{notation}{Notation and conventions}
\newtheorem*{contents}{Contents} 
\newtheorem{stepa}{Step}
\newtheorem{stepb}{Step}
\newtheorem{stepc}{Step}
\newtheorem{casea}{Case}

\newcommand{\lto}{\longrightarrow}

\begin{document}
\title{\bf Irreducibility of Hurwitz spaces}
\date{}
\author{Vassil Kanev}
\maketitle

\def\thefootnote{}
\footnote{Preprint N. 241, Dipartimento di Matematica, Universit\`{a} di Palermo.
}

\begin{abstract}
Graber, Harris and Starr proved, when $n\geq 2d$,  the irreducibility of 
the Hurwitz space
$\mathcal{H}^0_{d,n}(Y)$  which 
parametrizes degree $d$ coverings of  a smooth, projective curve $Y$ 
of genus $\geq 1$, simply branched in $n$ points, with full monodromy 
group $S_d$. We sharpen this result   and prove 
that $\mathcal{H}^0_{d,n}(Y)$ is irreducible if $n\geq \max\{2,2d-4\}$ and 
in the case of elliptic $Y$ if $n\geq \max\{2,2d-6\}$. We extend the result to coverings simply branched in all but one point of the discriminant. 
Fixing the
ramification multiplicities over the special point we prove that
the corresponding 
Hurwitz space is irreducible if the number of simply branched points 
is $\geq 2d-2$. We study also simply branched coverings with monodromy 
group $\neq S_d$ and when $n$ is large enough determine the corresponding connected components of $\mathcal{H}_{d,n}(Y)$. Our results are based on explicit calculation of the 
braid moves associated with the standard generators of the $n$-strand 
braid group of $Y$.
\end{abstract}

\section*{Introduction}
Let $Y$ be a smooth, connected, projective complex curve of genus 
$g\geq 0$. Let $\mathcal{H}_{d,n}(Y)$ be the Hurwitz  space 
which parametrizes degree $d$ coverings of $Y$   simply branched in $n$ points. A classical result of Hurwitz \cite{Hu} states that 
$\mathcal{H}_{d,n}(\mathbb{P}^1)$ is irreducible. More generally one 
can consider coverings of $\mathbb{P}^1$ which are simply branched in 
all but one point of the discriminant. Fixing the ramification multiplicities 
over the special point one obtains a corresponding Hurwitz space which turns 
out to be irreducible as well (see  \cite{Na},  \cite{Kl}, 
\cite{Mo}).

Coverings of curves of positive genus were studied by Graber, Harris 
and Starr in \cite{GHS}. They considered the Hurwitz space 
$\mathcal{H}^0_{d,n}(Y)$ which parametrizes coverings with full 
monodromy group $S_d$ and proved it is irreducible if $n\geq 2d$. 
Another result of this type was obtained by F. Vetro \cite{Ve}. 

In the present paper we sharpen the result of Graber, Harris and Starr 
and prove that $\mathcal{H}^0_{d,n}(Y)$ is irreducible if $n\geq 
2d-2$ (cf. Theorem~\ref{s4.45}). Our approach allows to extend the 
result to coverings which are simply branched in all but one point 
of the discriminant. Fixing the branching data of the special point, 
i.e. a partition 
$\underline{e}=\{e_1,e_2,\ldots,e_q\}$ where $e_1\geq e_2\geq \ldots 
\geq e_q\geq 1$ and $e_1+\cdots +e_q=d$ one obtains a Hurwitz space 
$\mathcal{H}^0_{d,n,\underline{e}}(Y)$ parametrizing coverings with 
full monodromy group, simply branched in $n$ points and ramified with 
multiplicities $e_1,\ldots ,e_q$ over one additional point. We prove 
$\mathcal{H}^0_{d,n,\underline{e}}(Y)$ is irreducible if $n\geq 2d-2$ 
(cf. Theorem~\ref{s4.49}). Coverings with $n<2d-2$ are more difficult 
to deal with. We prove the irreducibility of $\mathcal{H}^0_{d,n}(Y)$ 
when $n=2d-4,\; d\geq 3,\; g(Y)\geq 1$ (cf. Theorem~\ref{s5.64bis}) 
and the irreducibility of $\mathcal{H}^0_{d,n}(Y)$ when $n=2d-6,\; 
d\geq 4,\; g(Y)=1$ (cf. Theorem~\ref{s6.81}). The case of simply branched coverings 
with monodromy groups $\neq S_d$ can be reduced via \'{e}tale 
coverings to the case of coverings with full monodromy group. In 
Theorem~\ref{s4.55} and Theorem~\ref{s6.87} we prove that connected 
components of $\mathcal{H}_{d,n}(Y)$ not contained in 
$\mathcal{H}^0_{d,n}(Y)$ may exist only if $d$ is not prime and in 
this case if $n$ is sufficiently large such connected components 
correspond bijectively to the equivalence classes of \'{e}tale 
coverings $[\tilde{Y}\to Y]$ of degree $d_2$, where $d_2|d,\; d_2\neq 
1,d$. 

In proving these results we follow the standard approach for 
determining the connected components of $\mathcal{H}_{d,n}(Y)$. 
Associating to every equivalence class of coverings $[X\to Y]\in 
\mathcal{H}_{d,n}(Y)$ its 
discriminant locus  yields an \'{e}tale mapping 
$\mathcal{H}_{d,n}(Y)\to Y^{(n)}-\Delta$. Fixing a $D\in 
Y^{(n)}-\Delta$ the fiber over $D$ is identified via Riemann's 
existence theorem to the equivalence classes modulo inner 
automorphisms of ordered $(n+2g)$-tuples 
$(t_1,\ldots,t_n;\lambda_1,\mu_1,\ldots,\lambda_g,\mu_g)$ satisfying 
$t_{1}\cdots t_{n}=[\lambda_{1},\mu_{1}]\cdots [\lambda_{g},\mu_{g}]$ 
-- we call such $(n+2g)$-tuples \emph{Hurwitz systems} -- where 
$t_i,\lambda_k,\mu_k\in S_d$, \; $t_i$ are 
transpositions, and $t_i,\lambda_k,\mu_k$ with $1\leq i\leq n,\; 1\leq k\leq 
g$
generate a transitive subgroup of $S_d$. The problem is thus reduced 
to finding the orbits of the action of the $n$-strand braid group of 
$Y$, namely $\pi_1(Y^{(n)}-\Delta ,D)$, on the set of equivalence 
classes of Hurwitz systems. In fact it is more natural and it suffices 
to study the action of the braid group of the open Riemann surface 
$Y-\{b_0\}$, where $b_0$ is a fixed point. Birman found in 
\cite{Bi} a natural system of generators of these braid groups. Our 
results are based on Theorem~\ref{s2.27} where we calculate explicitly 
how these generators act on Hurwitz systems. 

We follow the 
key idea of \cite{GHS} which consists in applying a sequence of braid 
moves to a given Hurwitz system 
$(t_1,\ldots,t_n;\lambda_1,\mu_1,\ldots,\lambda_g,\mu_g)$ in order to replace 
it 
by a new one 
$(\tilde{t}_1,\ldots,\tilde{t}_n;\tilde{\lambda}_1,\tilde{\mu}_1,\ldots
,\tilde{\lambda}_g,\tilde{\mu}_g)$ in which as many as possible of the 
elements $\tilde{\lambda}_{k},\tilde{\mu}_{k}$ equal 1 and then 
reduce the sequence $(\tilde{t}_{1}\ldots ,\tilde{t}_{n})$ to a 
normal form using only elementary transformations of the Artin's braid 
group. An important tool in our arguments is  Main 
Lemma~\ref{s3.32} which states that if in a Hurwitz system 
$(t_1,\ldots,t_n;\lambda_1,\mu_1,\ldots,\lambda_g,\mu_g)$ one has 
$t_it_{i+1}=1$, then replacing the pair $(t_i,t_{i+1})$ by a pair 
$(h^{-1}t_ih,h^{-1}t_{i+1}h)$, where $h$ belongs to the subgroup generated by  
$t_j,\lambda_k,\mu_k$ with $1\leq j\leq n,\, j\neq i,i+1,\,1\leq k\leq g$, 
 one obtains a braid-equivalent Hurwitz system,  i.e. can be 
obtained from the initial one by a finite sequence of braid moves. This 
 is a generalization of Lemma~2.2 of \cite{GHS} proved there using the Kontsevich 
moduli space of stable maps. The  proof of  our Main 
Lemma~\ref{s3.32} is more elementary and uses only the explicit formulas 
for the braid moves of Theorem~\ref{s2.27}. We think  
Theorem~\ref{s2.27} and Main Lemma~\ref{s3.32} are results of independent 
interest, since both are valid for Hurwitz systems with values in an 
arbitrary (possibly infinite) group $G$. We profitted a lot from the 
paper of Mochizuki \cite{Mo} which was the starting point for us  in 
studying the problem of the irreducibility of Hurwitz spaces.
\begin{notation} 
We assume that the base field $k=\mathbb{C}$. 
Throughout the paper \linebreak
$\pi : X\to Y$ denotes a covering, i.e. a finite holomorphic mapping, of smooth projective curves.  
We do not assume that $X$ is connected, while $Y$ is always assumed connected.
Two coverings  $\pi_1 : X_1\to Y$ and $\pi_2 : X_2\to Y$ are called equivalent if there is an isomorphism $f:X_1\to X_2$ such that $\pi_1=\pi_2 \circ f$. A covering  $\pi : X\to Y$ of degree 
$d$ is called simple if both $X$ and $Y$ are irreducible  
 and for each $y\in Y$ one has $d-1\leq \#\; \pi^{-1}(y) \leq d$. We use mainly right actions and write $x^g$ for 
$x\in \Sigma$ and $g\in G$, where the group $G$ acts on the set $\Sigma$ on the right. In particular if $\Sigma=G$ we let $x^g=g^{-1}xg$. 
\end{notation}
\begin{contents}
1. Braid moves,\;
2. The main lemma,\;
3. The case $n\geq 2d-2$,
4. The case $n=2d-4,\quad g\geq 1$,\;
5. The case $n=2d-6,\: g=1$.\;
\end{contents}

\section{
Braid moves}
\label{s2}

\begin{block}\label{s2.11}
The equivalence classes of (possibly nonconnected) coverings of degree $d$ of 
a smooth, 
projective, irreducible curve $Y$ branched in $n>0$ points are 
parametrized by a smooth scheme $\mathcal{H}_d^{(n)}(Y)$ which is an 
\'{e}tale cover of $Y^{(n)}-\Delta$. Here $\Delta$ is the codimension 
one subscheme consisting of nonsimple divisors of degree $n$. The 
covering mapping associates to $[X\to Y]$ its branch locus $D$. Let 
$b_{0}\in Y$. We denote by $\mathcal{U}(b_0)$ 
 the open subset of $\mathcal{H}_d^{(n)}(Y)$ consisting of coverings which are 
unramified at 
$b_{0}$. Similarly we consider equivalence classes of pairs $[
\pi :X\to Y, \phi]$ where $[X\to Y]\in \mathcal{U}(b_0)$ and $\phi 
: \pi^{-1}(b_{0})\to \{1,\ldots ,d\}$ is a bijection. We denote by 
$\mathcal{H}_d^{(n)}(Y,b_0)$ the corresponding \'{e}tale cover of 
$\mathcal{U}(b_0)$. Let $D\in Y^{(n)}-\Delta,\; b_{0}\notin 
D$. Riemann's existence theorem (see e.g. \cite{Fu} p.544) 
establishes the following natural one-to-one correspondences
\renewcommand{\theenumi}{\Alph{enumi}}
\renewcommand{\theenumii}{\roman{enumii}}
\begin{enumerate}
\item For $\mathcal{H}_d^{(n)}(Y,b_0)$ between:
\begin{enumerate}
\item the fiber of $\mathcal{H}_d^{(n)}(Y,b_0)\to 
(Y- b_{0})^{(n)}- \Delta$ over $D$
\item the set of homomorphisms $m:\pi_1(Y- D,b_0)\to S_d$ 
which satisfy the condition that $m(\gamma)\neq 1$ for each small loop which 
circles a point of $D$.
\end{enumerate}
\item For $\mathcal{H}_d^{(n)}(Y)$ between: 
\begin{enumerate}
\item the fiber of 
$\mathcal{H}_d^{(n)}(Y)\to 
Y^{(n)}- \Delta$ over $D$
\item the set of equivalence classes modulo inner 
automorphisms of homomorphisms 
$m:\pi_1(Y - D,b_0)\to S_d$ satisfying the condition of A(ii), here $m\sim m'$ 
if there exists an $s\in S_d$ such that 
$m' = s^{-1}ms$.
\end{enumerate}
\end{enumerate}
We fix our preferences in this paper as follows. We consider the 
product in a fundamental group induced by product of arcs defined as 
\begin{equation}\label{es2.12}
\gamma_{1}*\gamma_{2}=
\begin{cases}
\gamma_{1}(2t)     & \text{for $0\leq t\leq 1/2$},\\
\gamma_{2}(2t-1)   & \text{for $1/2 \leq t\leq 1$}
\end{cases}
\end{equation}
We consider right action of $S_d$ on $\{1,\ldots ,d\}$, so a 
permutation $\sigma \in S_d$ transforms $i$ into $i^{\sigma}$ and for instance 
$(123)=(12)(13)$. Explicitly the correspondence in (A) is defined in 
the following way. Given a pair $\pi :X\to Y,\quad \phi :\pi^{-1}(b_0)\to 
\{1,\ldots,d\}$ let us number the elements of the fiber $\pi^{-1}(b_0)$ using 
$\phi$, thus $\pi^{-1}(b_0)=\{x_{1},\ldots , x_{d}\}$ where 
$x_{i}=\phi^{-1}(i)$. Then given a homotopy  class of loops $[\gamma]$ 
with $\gamma(0)=\gamma(1)=b_{0}$ the permutation $\sigma =m(\gamma)$ 
acts as follows. If lifting $\gamma$ one starts from $x_{i}$ and ends 
at $x_{j}$, then one lets $i^{\sigma }=j$. Let $s\in S_d$. Renumbering 
the points of the fiber $\pi^{-1}(b_0)$ by a bijective map $\psi 
:\pi^{-1}(b_0)\to \{1,\ldots ,d\}$ such that $\psi(x)=\phi(x)^{s}$ 
results in replacing the monodromy homomorphism $m$ by $s^{-1}ms$. 
Hence the one-to-one correspondence in (B) is obtained from that of (A) 
applying the forgetful mapping $[X\to Y,\; \phi]\mapsto [X\to Y]$.
\end{block}
\begin{block}\label{s2.12}
Suppose that $g(Y)\geq 1$. Let us fix the orientation of the real 2-manifold 
$Y$ considered as a complex manifold. Let $D=\{b_1,\ldots ,b_n\}\subset Y$. 
Let $b_0\notin D$. 
We describe a standard way of choosing generators of $\pi_1(Y-D,b_0)$. 
Cutting $Y$ along $2g$ simple closed arcs which begin at $b_0$ and do 
not contain any of $b_{i},\; i\geq 1$ one obtains a standard 
$4g$-polygon with sides 
$\alpha_1,\beta_1,\alpha_1^{-1},\beta_1^{-1},\ldots 
,\alpha_g,\beta_g,\alpha_g^{-1},\beta_g^{-1}$ which circle the polygon 
in counterclockwise direction. We consider a simple closed arc $L$ 
which begins at $b_{0}$,\quad $L-b_{0}$ is contained in the interior 
of the 4g-polygon and passes consecutively in counterclockwise 
direction through the points $\{b_1,\ldots ,b_n\}$. The closed arc $L$ 
divides the $4g$ polygon into two regions $R$ and $R'$ which stay on the 
left respectively on the right side of $L$ with respect to its counterclockwise 
orientation. We choose a simple arc $\ell_{1}$ which lies inside the 
region $R$ and connects $b_{0}$ and $b_{1}$. Then we choose a second 
simple arc $\ell_{2}$ inside $R$ which connects $b_{0}$ with $b_{2}$, 
has only $b_{0}$ as point in common with $\ell_{1}$, and lies on the 
left side of $\ell_{1}$. Continuing in this way we obtain an ordered 
$n$-tuple $(\ell_{1}\ldots ,\ell_{n})$ of simple arcs which do not 
meet outside $b_{0}$.  Let $\gamma_{i}$ be a closed path which 
begins at $b_{0}$, travels along $\ell_{i}$ to a point near $b_{i}$, 
makes a small counterclockwise loop around $b_{i}$, and returns to 
$b_{0}$ along $\ell_{i}$. We obtain a  $(n+2g)$-tuple of closed arcs 
$(\gamma_1,\ldots ,\gamma_n;\alpha_1,\beta_1,\ldots 
,\alpha_g,\beta_g)$ which we call a \emph{standard system of closed 
arcs}. The corresponding homotopy  classes yield a standard 
system of generators for $\pi_1(Y-D,b_0)$ which satisfy the only 
relation 
\begin{equation}\label{es2.13a}
\gamma_{1}\gamma_{2}\cdots \gamma_{n}\ \simeq \ 
[\alpha_{1},\beta_{1}]\cdots [\alpha_{g},\beta_{g}]
\end{equation}
Figure~\ref{fs2.13a} illustrates such a standard system. A reader who
prefers the clockwise orientation of closed arcs and ordering of the 
branch points from left to right may look at this and all subsequent 
figures from the other side of the sheet.
\begin{figure}[!ht]
\begin{center}
\includegraphics{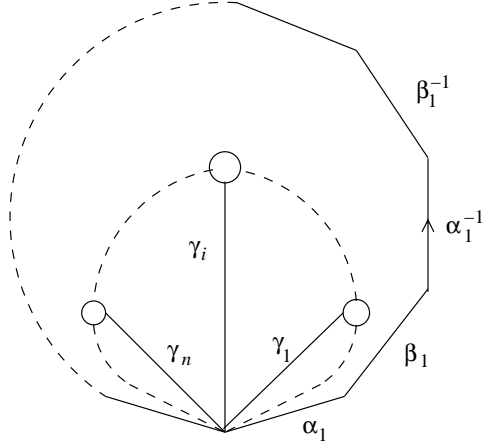}
\end{center}
\caption{Standard system of closed arcs}
\label{fs2.13a}
\end{figure}
Given a covering $\pi :X\to Y$ with discriminant $D$ and an 
isomorphism $\phi :\pi^{-1}(b_0)\to S_d$ one applies the monodromy 
homomorphism $m$ and obtains $t_{i}=m(\gamma_{i}),\; 
\lambda_{k}=m(\alpha_{k}),\; \mu_{k}=\lambda_{g+k}=m(\beta_{k})$.
\begin{dfn}\label{s2.14} An ordered sequence $(t_1,\ldots 
,t_n;\lambda_1,\mu_1,\ldots ,\lambda_g,\mu_g)$
of permutations in $S_d$ with $t_i\neq 1$ for $\forall i$ and 
satisfying the relation 
\[
t_1t_2\cdots t_n\ =\ [\lambda_{1},\mu_{1}]\cdots [\lambda_{g},\mu_{g}]
\]
is called a \emph{Hurwitz system}. We let $\lambda_{g+k}=\mu_{k}$. We call the 
subgroup $G\subseteq S_d$ 
generated by all $t_i,\lambda_{k},\mu_{k}$ the \emph{monodromy group} 
of the Hurwitz system.
\end{dfn}
By \ref{s2.11} given $D\in Y^{(n)}-\Delta $ and $b_{0}\notin D$, and 
fixing the closed arcs $(\gamma_1,\ldots 
,\gamma_n;\alpha_1,\beta_1,\ldots ,\alpha_g,\linebreak
\beta_g)$ as above, the 
fiber of $\mathcal{H}_d^{(n)}(Y,b_0)\to (Y-b_0)^{(n)}-\Delta $ over 
$D$ may be identified with the set of all Hurwitz systems, while the 
fiber of $\mathcal{H}_d^{(n)}(Y)\to Y^{(n)}-\Delta $ over $D$ may be 
identified with the set of equivalence classes 
$[t_1,\ldots ,t_n;\lambda_1,\mu_1,\ldots ,\lambda_g,\mu_g]$
of Hurwitz systems modulo inner automorphisms of $S_d$, where 
$(t_1,\ldots ,\mu_g)$ is equivalent to $(t'_{1},\ldots ,\mu'_{g})$ if 
there exists an $s\in S_d$ such that $t'_{i}=s^{-1}t_{i}s,\; 
\lambda'_{k}=s^{-1}\lambda_{k}s,\; \mu'_{k}=s^{-1}\mu_{k}s$ for 
$\forall i,k$.

An equivalent way of constructing a standard system of closed arcs is the 
following. One chooses first the $2g$ simple closed arcs 
$\alpha_1,\beta_1,\ldots,\alpha_g,\beta_g$. Then one chooses $n$ simple arcs 
which start at $b_0$, lie inside the $4g$-polygon, do not meet outside $b_0$, 
and have for end points the $n$ points of $D$. 
One enumerates these arcs 
according to the directions of departure in counterclockwise order.  
The 
obtained  $(n+2g)$-tuple $(\ell_1,\ldots,\ell_n;\alpha_1,\beta_1,\ldots,\alpha_g,\beta_g)$ is called an \emph{arc 
system} (cf. \cite{Lo} p.416).
One considers the induced ordering of the 
points of $D$. One can take for $R$ a star-like region which contains the 
union $\ell_1\cup\ldots \cup \ell_n$ and let $L=\partial R$. In this way one 
obtains all ingredients used to construct a standard system of closed arcs 
$(\gamma_1,\ldots ,\gamma_n;\alpha_1,\beta_1,\ldots 
,\alpha_g,\beta_g)$.
\end{block}
\begin{block}\label{s2.14a}
The connected components of $\mathcal{H}_d^{(n)}(Y)$ are in one-to-one 
correspondence with the orbits of the full $n$-strand braid group 
$\pi_1(Y^{(n)}-\Delta ,D)$ acting on the fiber  of the topological 
covering $\mathcal{H}_d^{(n)}(Y)\to Y^{(n)}-\Delta $ over $D$. Similar 
statement holds about $\mathcal{H}_d^{(n)}(Y,b_0)$ where one considers 
the braid group $\pi_1((Y-b_0)^{(n)}-\Delta ,D)$. The identification 
of these fibers with the Hurwitz systems reduces the problem of 
determining the connected components of $\mathcal{H}_d^{(n)}(Y,b_0)$ 
and $\mathcal{H}_d^{(n)}(Y)$ to calculating the action of the respective braid 
groups on Hurwitz systems and then finding the orbits.

Let $\gamma_{i},\alpha_{k},\beta_{k},\; 
1\leq i\leq n,\: 1\leq k\leq g$ be a standard system of closed arcs obtained from an arc system $\ell_{i},\alpha_{k},\beta_{k},\;
1\leq i\leq n,\: 1\leq k\leq g$
 as in \ref{s2.12}. Let $D^{u},\; 0\leq u\leq 1$ be a closed arc in $(Y-b_0)^{(n)}-\Delta 
$ with $D^{0}=D^{1}=D$. Suppose that starting from the given arc system one can extend the map $u\mapsto D^u$ to a homotopy of arc systems $\ell_{i}^u,\alpha_{k}^u,\beta_{k}^u,\; 1\leq i\leq n,\: 1\leq k\leq g$ based at $b_0$. This yields a corresponding homotopy of closed arcs $\gamma_{i}^u,\alpha_{k}^u,\beta_{k}^u
\; 1\leq i\leq n,\: 1\leq k\leq g$
 which form a standard system for each $u\in [0,1]$. Let $[X\to Y,\, \phi]\in 
\mathcal{H}_d^{(n)}(Y,b_0)$ have discriminant $D$ and let 
$(t_1,\ldots ,t_n;\lambda_1,\mu_1,\ldots 
,\lambda_g,\mu_g)$ be the corresponding Hurwitz system. Then the 
lifting of the closed arc $D^{u},\, 0\leq u\leq 1$, starting from a $[X\to Y,\: \phi]$, is defined by 
 by $m^{u}:\pi_1(Y-D^{u},b_0)\to S_d$ where (cf. 
\cite{Fu} p.545)
\[
m^{u}(\gamma_{i}^{u})=t_{i},\quad 
m^{u}(\alpha_{k}^{u})=\lambda_{k},\quad
m^{u}(\beta_{k}^{u})=\mu_{k}
\]
Letting $\gamma^{1}_{i}=\gamma'_{i},\; \alpha^{1}_{k}=\alpha'_{k},\; 
\beta^{1}_{k}=\beta'_{k}$ and $m^{1}=m'$ we obtain the end 
point of the lifting of $D^{u}$ is $[X'\to Y,\: \phi'] $ whose 
monodromy map $m':\pi_1(Y-D,b_0)\to S_{d}$ is defined by 
\begin{equation*}
m'(\gamma'_{i}) = t_{i},\quad m'(\alpha'_{k})=\lambda_{k},\quad 
m'(\beta'_{k})=\mu_{k}
\end{equation*}
Evaluating $m'$ at $\gamma_i,\alpha_k,\beta_k $ we obtain 
\begin{equation}\label{es2.16a}
m'(\gamma_{i})=t'_{i},\quad m'(\alpha_{k})=\lambda'_{k},\quad 
m'(\beta_{k})=\mu'_{k}.
\end{equation}
So $(t'_1,\ldots ,t'_n;\lambda'_1,\mu'_1,\ldots ,\lambda'_g,\mu'_g)$ 
is the Hurwitz system corresponding to $[X'\to Y,\: \phi'] $. Another 
approach is to consider 
\begin{equation}\label{es2.16c}
m(\gamma'_{i})=t''_{i},\quad 
m(\alpha'_{k})=\lambda''_{k},\quad m(\beta'_{k})=\mu''_{k}.
\end{equation}
 The Hurwitz 
system $(t''_1,\ldots ,t''_n;\lambda''_1,\mu''_1,\ldots ,\lambda''_g,\mu''_g)$
corresponds to a pair $[X''\to Y,\: \phi''] $ whose monodromy map 
$m'':\pi_1(Y-D,b_0)\to S_{d}$ satisfies
\begin{equation}\label{es2.16b}
m''(\gamma_{i})=t''_{i},\quad m''(\alpha_{k})=\lambda''_{k},\quad 
m''(\beta_{k}) = \mu''_{k}
\end{equation}
\end{block}
\begin{lem}\label{s2.16}
Let us lift the path $D^{1-u},\, 0\leq u\leq 1$ starting from $[X\to 
Y,\: \phi] $. Then the end point is $[X''\to Y,\: \phi''] $.
\end{lem}
\begin{proof}
In terms of $\gamma'_i,\alpha'_k,\beta'_k $ the monodromy map of 
$[X\to Y,\: \phi] $ is given by 
$m(\gamma'_{i})=t''_{i},\; 
m(\alpha'_{k})=\lambda''_{k},\; m(\beta'_{k})=\mu''_{k}$. Let us 
consider the homotopy of closed paths 
$\gamma_i^{1-u},\alpha_k^{1-u},\beta_k^{1-u}$, $0\leq u\leq 1 $ and 
the path in $\mathcal{H}_d^{(n)}(Y,b_0)$ given by 
$n^{u}:\pi_1(Y-D^{1-u},b_0)\to S_{d}$ where 
$n^{u}(\gamma^{1-u}_{i})=t_{i}'',\, 
n^{u}(\alpha^{1-u}_{k})=\lambda''_{k},\, 
n^{u}(\beta^{1-u}_{k})=\mu''_{k}$. Then $n^{0}=m,\; n^{1}=m''$ (cf. 
Eq.\eqref{es2.16b})
\end{proof}
\begin{dfn}\label{s2.17}
Given a closed arc $D^{u},\, 0\leq u\leq 1$ in the configuration space 
$(Y-b_0)^{(n)}-\Delta $ the transformation of Hurwitz systems 
$(t_1,\ldots ,t_n;\lambda_1,\mu_1,\ldots ,\lambda_g,\mu_g)\mapsto 
(t'_1,\ldots ,t'_n;\lambda'_1,\mu'_1,\ldots ,\lambda'_g,\linebreak
\mu'_g)$ given 
by Eq.\eqref{es2.16a} is called a \emph{braid move} of the first type. 
The transformation $(t_1,\ldots ,t_n;\lambda_1,\linebreak
\mu_1,\ldots ,\lambda_g,\mu_g)\mapsto (t''_1,\ldots ,t''_n;\lambda''_1,\mu''_1,\ldots 
,\lambda''_g,\mu''_g)$ given by Eq.\eqref{es2.16b} is called a braid move 
of the second type. 
\end{dfn}
The braid moves of the first and of the second 
type are inverse to each other according to Lemma~\ref{s2.16}.
It is evident that the braid moves of both types commute with inner 
automorphisms of $S_{d}$. So the braid moves are well-defined on 
equivalence classes of Hurwitz systems (cf.  \ref{s2.12}).
\begin{block}\label{s2.17a}
There is a convenient system of generators of 
$\pi_1((Y-b_0)^{(n)}-\Delta ,D)$. We include here some 
material borrowed from \cite{Bi} and \cite{Sc} for the sake of 
convenience of the reader and since our choices differ slightly from 
theirs. Consider the Galois covering 
$p:(Y-b_0)^{n}\to (Y-b_0)^{(n)}$ with Galois group $S_{n}$. 
Restricting to the 
complement of $\Delta$  one obtains an unramified 
Galois covering $p:(Y-b_0)^{n}-p^{-1}(\Delta)\to 
(Y-b_0)^{(n)}-\Delta$. Following the notation of \cite{FN} if 
$Q_{1}=\{b_{0}\}$ one denotes $(Y-b_0)^{n}-p^{-1}(\Delta)$ by 
$F_{1,n}Y$. Let $D=\{b_1,\ldots,b_n\},\; 
\tilde{D}=(b_1,\ldots,b_n)$.
One has an exact sequence
\begin{equation}\label{es2.17a}
1\lto \pi_1(F_{1,n}Y,\tilde{D})\lto \pi_1((Y-b_0)^{(n)}-\Delta 
,D)\lto S_{n}\lto 1
\end{equation}
One determines first a system of generators of the pure braid group
$\pi_1(F_{1,n}Y,\tilde{D})$
 as follows. Consider the closed arcs in $F_{1,n}Y$ defined by 
$(b_1,\ldots,b_{i-1},r_{ik}(t),b_{i+1},\ldots ,b_n)$ and 
$(b_1,\ldots,b_{i-1},t_{ik}(t),b_{i+1},\ldots ,b_n)$, with $t\in [0,1]$, 
where $r_{ik}$ and $t_{ik}$ are the closed simple arcs based at 
$b_{i}$ and pictured on Figure~\ref{fs2.17a}   by a continuous 
line and by a dotted line respectively.
\begin{figure}[!ht]
\begin{center}
\includegraphics{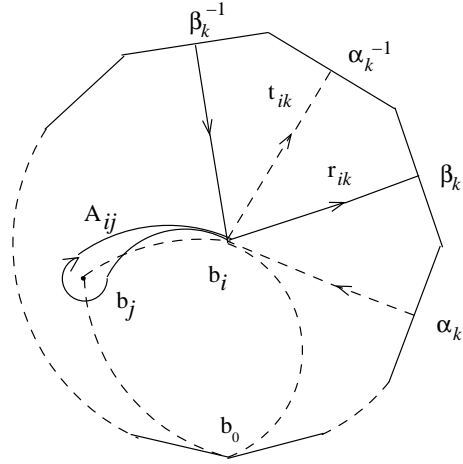}
\end{center}
\caption{Generators of the pure braid group $\pi_1(F_{1,n},\tilde{D})$}
\label{fs2.17a}
\end{figure}
We denote the corresponding homotopy  classes by $\rho_{ik} , 
\tau_{ik} \in \pi_1(F_{1,n}(Y),\tilde{D})$. Informally $\rho_{ik} $ 
corresponds to a loop of the $i$-th point along $\alpha_{k}$ and 
$\tau_{ik} $ corresponds to a loop of the $i$-th point along 
$\beta_{k}$. Let us denote by $A_{ij},\, i<j$ the element of 
$\pi_1(F_{1,n}Y,\tilde{D})$ represented by a closed simple arc in 
$Y^{n}$ which leaves fixed $b_{k}$ for $k\neq i$ and moves $b_{i}$ 
along the arc pictured on Figure~\ref{fs2.17a}. For every $i$ and 
every $j>i$ the pictured loop is chosen so as to stay on the left of 
all arcs used to construct $\rho_{ik} $ and $\tau_{ik} $ for 
$k=1,\ldots ,g$. 
\begin{claim}
The pure braid group $\pi_1(F_{1,n}Y,\tilde{D})$ is generated by 
$\rho_{ik} ,\tau_{ik} $ and $A_{ij}$ where $i,j=1,\ldots,n,\quad 
i<j$ and $k=1,\ldots ,g$.
\end{claim}
\begin{proof}
This is proved by induction on $n$. When $n=1$ the claim is obvious. 
Let $n\geq 2$. Consider the fibration $F_{1,n}Y\to F_{1,n-1}Y$ defined 
by $(y_1,y_2,\ldots,y_n)\to (y_2,\ldots ,y_n)$. One has an exact 
sequence
\[
1 \to \pi_1(Y-\{b_0,b_2,\ldots,b_n\},b_1)\to 
\pi_1(F_{1,n}Y,(b_1,\ldots,b_n))\to 
\pi_1(F_{1,n-1}Y,(b_2,\ldots,b_n))\to 1
\]
since $\pi_2(F_{1,n-1}Y)=1$ by \cite{FN}, Corollary~2.2.
The elements $\rho_{ik} ,\tau_{ik} ,A_{ij}$ with $i\geq 2$ map to 
the corresponding elements in $\pi_1(F_{1,n-1}Y,(b_2,\ldots,b_n))$. 
The group $\pi_1(Y-\{b_0,b_2,\ldots,b_n\},b_1)$ is freely generated by 
elements which map into $\rho_{1,k} ,\tau_{1,k} ,A_{1,j}$ where $j\geq 
2,\quad k=1,\ldots ,g$. This shows the claim.
\end{proof}
Following \cite{Fu} p.547 and using the notation of \ref{s2.12} let us 
denote by $R_{i}$ the simply connected region of $R$ enclosed by the 
arcs $\ell_{i}$ and $\ell_{i+1}$ and the arc of $L$ from $b_i$ to 
$b_{i+1}$. For every $i=1,\ldots ,n-1$ choose simple arcs 
$\eta_{i}:[0,1]\to Y$ going from $b_{i}$ to $b_{i+1}$ in $R_{i}$ and 
$\eta'_{i}:[0,1]\to Y$ going from $b_{i+1}$ to $b_{i}$ in $R'$. Let 
$s_{i}:[0,1]\to Y^{(n)}-\Delta$ be the closed arc 
\[
s_{i}(t)\ =\ 
\{b_1,\ldots,b_{i-1},\eta_{i}(t),\eta'_{i}(t),b_{i+2},\ldots,b_n\}
\]
The homotopy  class of $s_{i}$ is denoted by $\sigma_{i}$. We may 
consider $\pi_1(F_{1,n}(Y),\tilde{D})$ as embedded in 
$\pi_1((Y-b_0)^{(n)}-\Delta ,D)$. The following relations are easy to 
verify (cf. \cite{FB} p.249)
\begin{equation}\label{es2.19}
\begin{split}
&\rho_{i+1,k}\ =\ \sigma_{i}\rho_{ik} \sigma_{i}^{-1},\quad 
\tau_{i+1,k} \
=\ \sigma_{i}\tau_{ik} \sigma_{i}^{-1}\\
&A_{ij}\ =\ \sigma _{i}^{-1}\cdots \sigma _{j-2}^{-1}\sigma 
_{j-1}^{2}\sigma _{j-2}\cdots \sigma _{i}
\end{split}
\end{equation}
\end{block}
\begin{pro}[Birman]\label{s2.19}
Let $a,b\in \mathbb{Z},\quad 1\leq a,b\leq n$. The braid group \linebreak
$\pi_1((Y-b_0)^{(n)}-\Delta ,D)$ is generated by 
$\sigma _{j},\rho_{ak},\tau_{bk}$ where $1\leq j\leq n-1,\quad 1\leq 
k\leq g$. The corresponding homotopy  classes generate 
$\pi_1(Y^{(n)}-\Delta ,D)$ as well.
\end{pro}
\begin{proof}
Let us consider the exact sequence Eq.\eqref{es2.17a}. The braids $\sigma 
_{j}$ map to the transpositions $(j\: j+1),\quad j=1,\ldots,n-1$ which 
generate $S_{n}$. So $\sigma _{j},\rho_{ik},\tau_{ik},A_{ij}$ with 
$i<j$ generate $\pi_1((Y-b_0)^{(n)}-\Delta ,D)$ according to the claim 
proved in  \ref{s2.17a}. The relations Eq.\eqref{es2.19} show that the 
generators may be reduced as stated in the theorem. The last statement 
follows from the surjection $\pi_1((Y-b_0)^{(n)}-\Delta 
,D)\twoheadrightarrow \pi_1(Y^{(n)}-\Delta ,D)$.
\end{proof}
We described in \ref{s2.17a} closed arcs in $(Y-b_0)^{(n)}-\Delta $ 
based at $D=\{b_1,\ldots,b_n\}$ whose homotopy classes form a system 
of generators $\sigma _{j},\rho_{ik},\tau_{ik}$ for 
$\pi_1((Y-b_0)^{(n)}-\Delta ,D)$. Our aim now is for each of these to 
construct a homotopy of the standard system of closed 
arcs $\gamma_{i}^{u},\alpha_{k}^{u},\beta_{k}^{u}$ as in \ref{s2.14a}. 
This will permit us to calculate eventually the corresponding braid 
moves of the Hurwitz systems. The calculation of the braid moves 
$\sigma _{j},\: j=1,\ldots,n-1$ is due to Hurwitz (cf. \cite{Hu} or 
e.g. \cite{Vo}, Theorem~10.3). We define closed arcs 
$\delta_{k}:[0,1]\to Y,\: \delta_k(0)=\delta_k(1)=b_{0},\:
k=0,1,\ldots,g$ as follows. We let $\delta_{0}(t)=b_{0},\: \forall 
t\in[0,1]$. We connect the initial vertex of $\alpha_{1}$ with the end 
vertex of $\beta_{1}^{-1}$ in the $4g$-polygon of Figure~\ref{fs2.13a} by a 
simple arc which belongs to the region $R'$ on the right of $L$ (cf. 
\ref{s2.12}). This yields $\delta_{1}$. We connect the initial vertex 
of $\alpha_{1}$ with the end vertex of $\beta_{2}^{-1}$ by a simple 
arc which belongs to the region on the right of $L$ and on the left of 
$\delta_{1}$. We denote the corresponding closed arc of $Y$ by 
$\delta_{2}$. Continuing in this way we obtain 
$\delta_{1},\delta_{2},\ldots,\delta_{g}$ (see Figure~\ref{fs2.22}). 
Clearly $\delta_{k}\simeq [\alpha_{1},\beta_{1}]\cdots 
[\alpha_{k},\beta_{k}]$ in $Y-D$ and  $\delta_{g}\simeq 
\gamma_{1}\cdots \gamma_{n}$ according to Eq.\eqref{es2.13a}.
\begin{thm}\label{s2.20}
Let $Y$ be a compact, closed Riemann surface of genus $g(Y)\geq 1$. 
Let $b_{0}\in Y$ and let $\gamma_{i},\alpha_{k},\beta_{k}$ with $1\leq 
i\leq n,\: 1\leq k\leq g$ be a standard system of closed arcs as in 
\ref{s2.12}. Let $\delta_{k},\: k=0,1,\ldots ,g$ be the closed arcs 
defined above. For each of the closed arcs in $(Y-b_0)^{(n)}-\Delta $ 
constructed in \ref{s2.17a} and representing $\sigma 
_{j},\rho_{ik},\tau_{ik}$ there is a homotopy
$\gamma_{i}^{u},\alpha_{k}^{u},\beta_{k}^{u},\: u\in [0,1]$ of the 
standard system of closed arcs such that the end system 
$\gamma'_{i},\alpha'_{k},\beta'_{k},\: 1\leq i\leq n,\: 1\leq 
k\leq g$ is homotopic to:

\smallskip
a.\quad for $\sigma_{j}$ where $1\leq j\leq n-1$
\begin{gather*}
\gamma'_{i}\simeq \gamma_{i}\quad \text{for $\forall i\neq j,\, j+1$}, 
\quad \alpha'_{k}\simeq \alpha_{k},\; \beta'_{k}\simeq \beta_{k}\quad 
\text{for $\forall k$}\\
\gamma'_{j}\simeq \gamma_{j+1}, \quad \gamma'_{j+1}\simeq 
\gamma^{-1}_{j+1}\gamma_{j}\gamma_{j+1}
\end{gather*}

b.\quad for $\rho_{ik}$ where $1\leq i\leq n,\: 1\leq k\leq g$
\begin{gather*}
\gamma'_{j}\simeq \gamma_{j}\quad \text{for $\forall j\neq i$},\quad
\alpha'_{\ell}\simeq \alpha_{\ell}\quad \text{for $\forall \ell$}, 
\quad \beta'_{\ell}\simeq \beta_{\ell}\quad \text{for $\forall \ell 
\neq k$}\\
\gamma'_{i}\simeq \mu_{ik}\gamma_{i}\mu_{ik}^{-1},\quad 
\beta'_{k}\simeq (\zeta_{ik}\gamma_{i}\zeta_{ik}^{-1})\beta_{k},\quad 
\text{where}\\
\mu_{ik}\simeq (\gamma_{1}\cdots 
\gamma_{i-1})^{-1}\delta_{k-1}\alpha_{k}(\delta_{k}^{-1}\delta_{g})
(\gamma_{i+1}\cdots \gamma_{n})^{-1}, \quad 
\zeta_{ik}\simeq (\delta_{k}^{-1}\delta_{g})(\gamma_{i+1}\cdots 
\gamma_{n})^{-1}
\end{gather*}

c.\quad for $\tau_{ik}$ where $1\leq i\leq n,\: 1\leq k\leq g$
\begin{gather*}
\gamma'_{j}\simeq \gamma_{j}\quad \text{for $\forall j\neq i$},\quad
\alpha'_{\ell}\simeq \alpha_{\ell}\quad \text{for $\forall \ell \neq k$}, 
\quad \beta'_{\ell}\simeq \beta_{\ell}\quad \text{for $\forall \ell $}\\
\gamma'_{i}\simeq \nu_{ik}\gamma_{i}\nu_{ik}^{-1},\quad 
\alpha'_{k}\simeq (\xi_{ik}
\gamma_{i}^{-1}\xi_{ik}^{-1})\alpha_{k},\quad 
\text{where}\\
\nu_{ik}\simeq \gamma_{i+1}\cdots 
\gamma_{n}(\delta_{k}^{-1}\delta_{g})^{-1}\beta_k 
\delta_{k-1}^{-1}\gamma_{1}\cdots 
\gamma_{i-1}, \quad \xi_{ik}\simeq \delta_{k-1}^{-1}\gamma_{1}\cdots 
\gamma_{i-1}
\end{gather*}
For each of the inverse closed arcs corresponding to 
$\sigma_{j}^{-1},\rho_{ik}^{-1},\tau_{ik}^{-1}$ there is a 
corresponding homotopy of the standard system 
$\gamma_{i},\alpha_{k},\beta_{k}$, $1\leq i\leq n,\: 1\leq k\leq g$
such that the end system $\gamma''_{i},\alpha''_{k},\beta''_{k}$,
$1\leq i\leq n,\: 1\leq k\leq g$
is homotopic to:

\smallskip
d.\quad for $\sigma_{j}^{-1}$ where $1\leq j\leq n-1$
\begin{gather*}
\gamma''_{i}\simeq \gamma_{i}\quad \text{for $\forall i\neq j,\, j+1$}, 
\quad \alpha''_{k}\simeq \alpha_{k},\; \beta''_{k}\simeq \beta_{k}\quad 
\text{for $\forall k$}\\
\gamma''_{j}\simeq \gamma_{j}\gamma_{j+1}\gamma^{-1}_{j}, \quad 
\gamma''_{j+1}\simeq \gamma_{j} 
\end{gather*}

e.\quad for $\rho_{ik}^{-1}$ where $1\leq i\leq n,\: 1\leq k\leq g$
\begin{gather*}
\gamma''_{j}\simeq \gamma_{j}\quad \text{for $\forall j\neq i$},\quad
\alpha''_{\ell}\simeq \alpha_{\ell}\quad \text{for $\forall \ell$}, 
\quad \beta''_{\ell}\simeq \beta_{\ell}\quad \text{for $\forall \ell 
\neq k$}\\
\gamma''_{i}\simeq \mu_{ik}^{-1}\gamma_{i}\mu_{ik},\quad 
\beta''_{k}\simeq 
(\tilde{\zeta}_{ik}^{-1}\gamma_{i}^{-1}\tilde{\zeta}_{ik})\beta_{k},\quad 
\text{where}\\
\text{$\mu_{ik}$ is as in (b) and}\quad 
\tilde{\zeta}_{ik}\simeq (\gamma_{1}\cdots 
\gamma_{i-1})^{-1}\delta_{k-1}\alpha_{k}
\end{gather*}

f.\quad for $\tau_{ik}^{-1}$ where $1\leq i\leq n,\: 1\leq k\leq g$
\begin{gather*}
\gamma''_{j}\simeq \gamma_{j}\quad \text{for $\forall j\neq i$},\quad
\alpha''_{\ell}\simeq \alpha_{\ell}\quad \text{for $\forall \ell \neq k$}, 
\quad \beta''_{\ell}\simeq \beta_{\ell}\quad \text{for $\forall \ell $}\\
\gamma''_{i}\simeq \nu_{ik}^{-1}\gamma_{i}\nu_{ik},\quad 
\alpha''_{k}\simeq (\tilde{\xi}_{ik}^{-1}
\gamma_{i}\tilde{\xi}_{ik})\alpha_{k},\quad 
\text{where}\\
\text{$\nu_{ik}$ is as in (c) and}\quad
\tilde{\xi}_{ik}\simeq (\gamma_{i+1}\cdots 
\gamma_{n})(\delta_{k}^{-1}\delta_{g})^{-1}\beta_{k}
\end{gather*}
\end{thm}
\begin{proof}
(a)\quad Here one moves only $\gamma_{j}$ and $\gamma_{j+1}$. Clearly 
$\gamma'_{j}\simeq \gamma_{j+1}$ and $\gamma_{j}\gamma_{j+1}\simeq 
\gamma'_{j}\gamma'_{j+1}$, so $\gamma'_{j+1}\simeq 
\gamma^{-1}_{j+1}\gamma_{j}\gamma_{j+1}$.

\smallskip
(b)\quad The effect of the homotopy of 
$\gamma_{i},\alpha_{k},\beta_{k}$ along $\rho_{ik}$ is pictured on 
Figure~\ref{fs2.22}.
\begin{figure}[!ht]
\begin{center}
\scalebox{.90}{\includegraphics{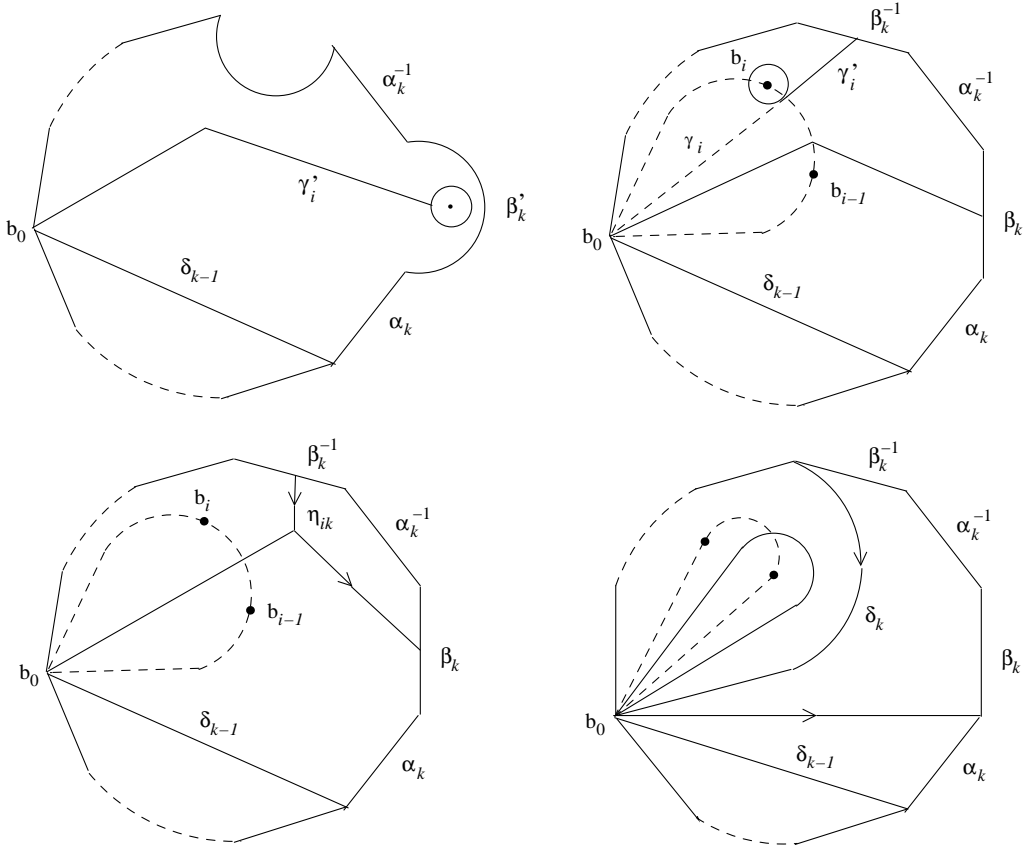}}
\end{center}
\caption{Homotopy along $\rho_{ik}$}
\label{fs2.22}
\end{figure}
One moves the point $b_{i}$ along the arc $r_{ik}$ pictured on 
Figure~\ref{fs2.17a} and together with it deforms the closed arc 
$\gamma_{i}$. At the moment $\gamma_{i}^{u}$ reaches the side 
$\beta_{k}$ one deforms also $\beta_{k}$ in order that the condition 
$\gamma_{i}^{u}$ and $\beta_{k}^{u}$ have no points in common except 
$b_{0}$ remains valid. None of $\gamma_{j}$ for $j\neq i$, or 
$\alpha_{\ell}$ for $\forall \ell$, or $\beta_{\ell}$ for $\forall 
\ell \neq k$ changes in this homotopy. The effect of cutting $Y$ 
along the closed arcs 
$\alpha_{1},\beta_{1},\ldots,\alpha_{k},\beta'_{k},\ldots,\alpha_{g},
\beta_{g}$ is the same as to cut a region containing $b_{i}$ from 
the original $4g$-polygon and glue it along the side $\beta_{i}$ as 
described in Figure~\ref{fs2.22}~(NW) (NW=Northwest). We wish to express $\gamma'_{i}$ in 
terms of the standard system $\gamma_{j},\alpha_{\ell},\beta_{\ell}$, 
$1\leq j\leq n, 1\leq \ell\leq g$. The closed arc $\gamma'_{i}$ is 
pictured in the original $4g$-polygon on Figure~\ref{fs2.22}~(NE).
It  is  clear  that  $\gamma'_{i}\simeq  \eta_{ik}\gamma_{i}\eta_{ik}^{-1}$ 
where    $\eta_{ik}$    is    the    closed    arc     pictured     on 
Figure~\ref{fs2.22}~(SW).   That    $\eta_{ik}$    is    homotopic    to 
$(\gamma_{1}\cdots 
\gamma_{i-1})^{-1}(\delta_{k-1}\alpha_{k})\delta_{k}^{-1}(\gamma_{1}
\cdots  \gamma_{i-1})$   is   evident   from   Figure~\ref{fs2.22}~(SE). 
Furthermore $\eta_{ik}\gamma_{i}\eta_{ik}^{-1}\simeq 
(\eta_{ik}\gamma_{i})\gamma_{i}(\eta_{ik}\gamma_{i})^{-1}$. Using  the 
relation  $\gamma_{1}\cdots  \gamma_{n}\simeq  \delta_{g}$  we  obtain 
$\eta_{ik}\gamma_{i}\simeq
(\gamma_{1}\cdots 
\gamma_{i-1})^{-1}\delta_{k-1}\alpha_{k}(\delta_{k}^{-1}\delta_{g})
(\gamma_{i+1}\cdots \gamma_{n})^{-1}$.
This proves the formula for $\gamma'_{i}$ of Part~(b). The calculation of 
$\beta'_{k}$ is similar. The closed arc $\beta'_{k}$ is pictured in 
the original $4g$-polygon on Figure~\ref{fs2.25}~(NW).
\begin{figure}[!ht]
\begin{center}
\includegraphics{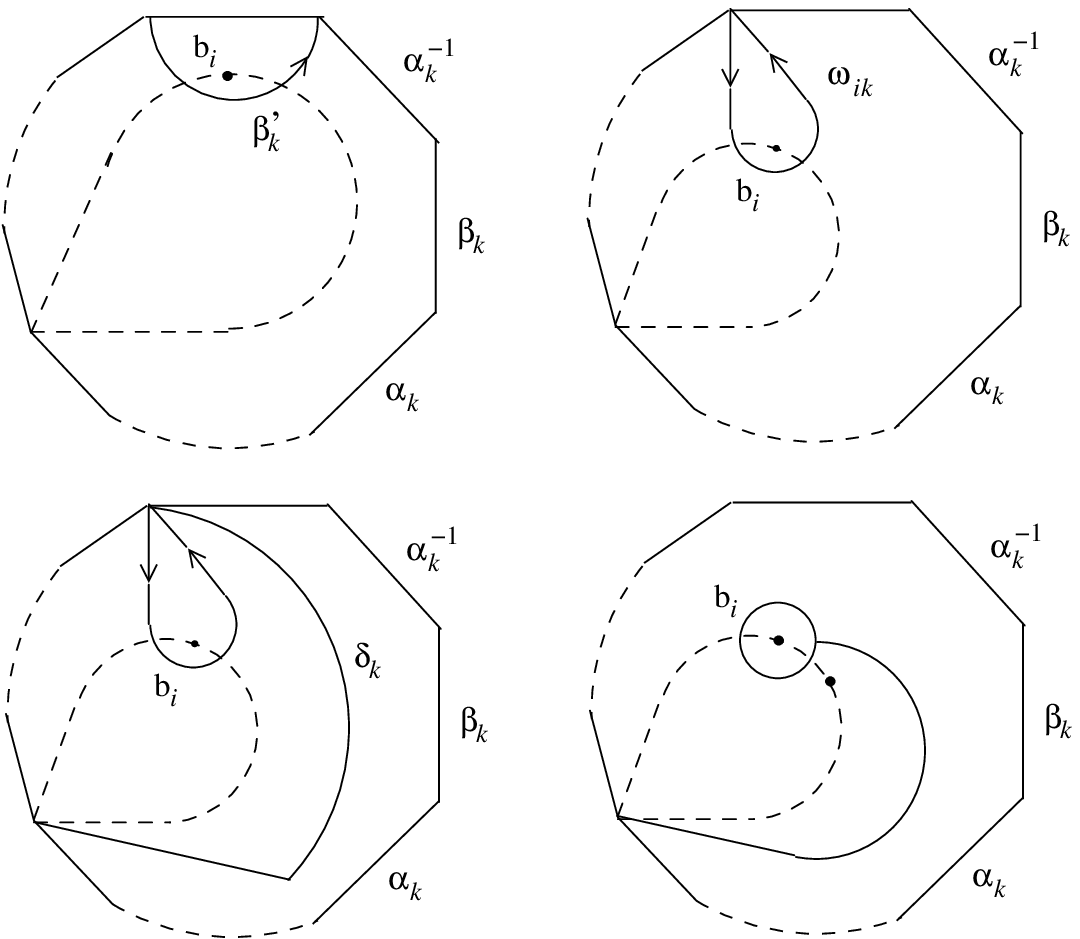}
\end{center}
\caption{{}}
\label{fs2.25}
\end{figure}
It is homotopic 
to $\omega_{ik}\cdot \beta_{k}$ where $\omega_{ik}$ is the closed arc 
based at $b_{0}$ pictured on Figure~\ref{fs2.25}~(NE). We then consider 
$\delta_{k}\omega_{ik}\delta_{k}^{-1}$ (see Figure~\ref{fs2.25}~(SW)). 
The latter is homotopic to $(\gamma_{1}\cdots 
\gamma_{i-1})\gamma_{i}(\gamma_{1}\cdots \gamma_{i-1})^{-1}$. We thus 
obtain  $\omega_{ik}\simeq (\delta_{k}^{-1}\gamma_{1}\cdots 
\gamma_{i-1})\gamma_{i}(\gamma_{1}\cdots \gamma_{1})^{-1}\delta_{k}$. We then 
have 
\begin{align*}
\omega_{ik}
&\simeq \delta_{k}^{-1}(\gamma_{1}\cdots 
\gamma_{i-1}\gamma_{i})\gamma_{i}(\gamma_{1}\cdots 
\gamma_{i-1}\gamma_{i})^{-1}\delta_{k}\\
&\simeq [\delta_{k}^{-1}\delta_{g}(\gamma_{i+1}\cdots 
\gamma_{n})^{-1}]\gamma_{i}[\gamma_{i+1}\cdots \gamma_{n}
\delta_{g}^{-1}\delta_{k}]
\end{align*}
since $\gamma_{1}\cdots \gamma_{n}\simeq \delta_{g}$. This proves the 
second formula of Part~(b).

\smallskip
(c)\quad The arguments here are very similar to those of 
Part~(b). One 
deforms $\gamma_{i}$ and $\alpha_{k}$ with $b_{i}^{u}$ moving along 
the arc $t_{ik}$ (see Figure~\ref{fs2.17a}). When $b_{i}^{u}$ returns 
to $b_{i}$ one obtains closed arcs $\gamma_{i}'$ and $\alpha_{k}'$ 
for which $\gamma_{i}'\simeq \theta_{ik}\gamma_{i}\theta_{ik}^{-1}$, 
$\alpha_{k}'=\varepsilon_{ik}\alpha_{k}$ where $\theta_{ik}$ and $\varepsilon_{ik}$ 
are represented by arcs in the original $4g$-polygon pictured on 
Figure~\ref{fs2.26}. 
\begin{figure}[!ht]
\begin{center}
\includegraphics{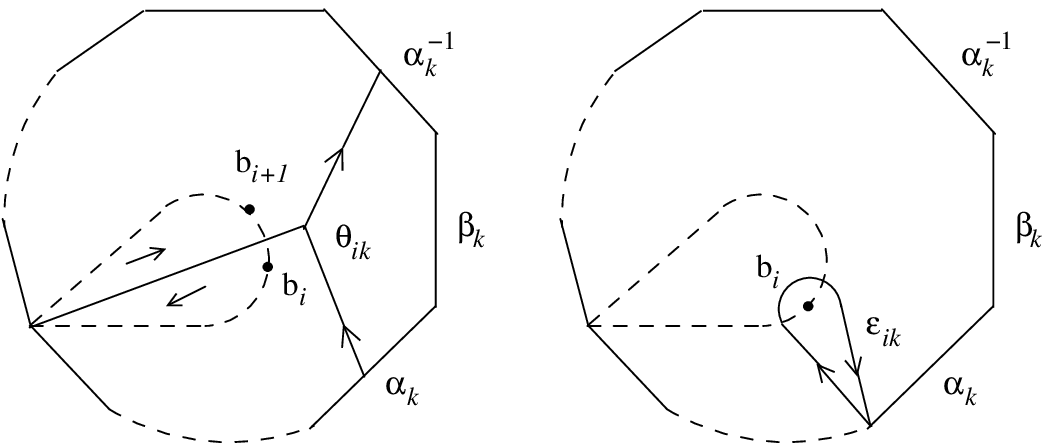}
\end{center}
\caption{{}}
\label{fs2.26}
\end{figure}
We calculate $\theta_{ik}$ in the following way. We consider  it as a 
product of four arcs according to the picture. We deform $\theta_{ik}$ 
in such a way that the first arc becomes a closed simple arc encircling 
$\{b_1,\ldots,b_i\}$ in clockwise direction, then the second arc goes 
from the initial point of $\alpha_{1}$ to the end point of 
$\alpha_{k}^{-1}$, the third arc goes from the initial point of 
$\alpha_{k}$ to the initial point of $\alpha_{1}$ and the fourth arc 
equals the first one with the opposite orientation. We have accordingly
\[
\theta_{ik} \simeq (\gamma_{1}\cdots 
\gamma_{i})^{-1}(\delta_{k}\beta_{k})\delta_{k-1}^{-1}(\gamma_{1}\cdots 
\gamma_{i})
\]
Conjugating $\gamma_{i}$ by $\theta_{ik}$ we may cancel the last 
factor $\gamma_{i}$ from $\theta_{ik}$ and replace $(\gamma_{1}\cdots 
\gamma_{i})^{-1}$ by $(\gamma_{i+1}\cdots \gamma_{n})\delta_{g}^{-1}$. 
We thus obtain $\gamma_{i}'=\nu_{ik}\gamma_{i}\nu_{ik}^{-1}$ where 
\[
\nu_{ik}\ \simeq \ (\gamma_{i+1}\cdots 
\gamma_{n})(\delta_{k}^{-1}\delta_{g})^{-1}\beta_{k}\delta_{k-1}^{-1}(
\gamma_{1}\cdots \gamma_{i-1})
\]
Finally, $\delta_{k-1}\varepsilon_{ik}\delta_{k-1}^{-1}\simeq 
(\gamma_{1}\cdots \gamma_{i-1})\gamma_{i}^{-1}(\gamma_{1}\cdots 
\gamma_{i-1})^{-1}$.
This proves the last formula of Part~(c).

In order to obtain the formulas of (d), (e) and (f) from those of (a), 
(b) and (c) we notice that if we apply in each case to 
$(\gamma_{1}'',\ldots ,\gamma_{n}'';\alpha_{1}'',\ldots,\beta_{g}'')$ 
the braid moves $\sigma_{j},\rho_{ik},\tau_{ik}$ respectively  we 
obtain 
$(\gamma_{1},\ldots ,\gamma_{n};\alpha_{1},\ldots,\beta_{g})$. For 
example, in order to verify (d) we have $\gamma_{j}\simeq 
\gamma_{j+1}'',\: \gamma_{j+1}\simeq 
(\gamma_{j+1}'')^{-1}\gamma_{j}''\gamma_{j+1}''$. Therefore 
$\gamma_{j}''\simeq \gamma_{j}\gamma_{j+1}\gamma_{j}^{-1}$. In  
Case (e) applying $\rho_{ik}$ to $(\gamma_{1}'',\ldots, \beta_{g}'')$ 
we obtain $\gamma_{j}\simeq \gamma_{j}''$ for $\forall j\neq i,\: 
\alpha_{\ell}\simeq \alpha_{\ell}''$ for $\forall \ell$, 
$\beta_{\ell}\simeq \beta_{\ell}''$ for $\forall \ell \neq k$. 
Furthermore $\gamma_{i}\simeq \mu_{ik}''\gamma_{i}''(\mu_{ik}'')^{-1}$ 
where $\mu_{ik}''$ is expressed by 
$\gamma_{j}'',\alpha_{\ell}'',\beta_{\ell}''$ as in (b). Since neither 
$\gamma_{i}''$ nor $\beta_{k}''$ enter in this expression we may replace 
$\gamma_{j}'',\alpha_{\ell}'',\beta_{\ell}''$ by 
$\gamma_{j},\alpha_{\ell},\beta_{\ell}$ and we obtain 
$\mu_{ik}''=\mu_{ik}$. Thus $\gamma_{i}''\simeq 
\mu_{ik}^{-1}\gamma_{i}\mu_{ik}$. Similarly we have 
$\beta_{k}=(\zeta_{ik}\gamma_{i}''\zeta_{k}^{-1})\beta_{k}''$ where 
$\zeta_{ik}$ is as in (b). Replacing $\gamma_{i}''$ by 
$\mu_{ik}^{-1}\gamma_{i}\mu_{ik}$ and canceling we obtain 
$\beta_{k}''=(\tilde{\zeta}_{ik}^{-1}\gamma_{i}^{-1}\tilde{\zeta}_{ik}
)\beta_{k}$ where $\tilde{\zeta}_{ik}$ is as in (e). In a similar 
manner one deduces (f) from (c).
\end{proof}
Recall from \ref{s2.11} and \ref{s2.12} that given $D\in 
Y^{(n)}-\Delta ,\; b_{0}\in Y-D$ and fixing a standard system of 
closed arcs 
$\gamma_{1},\ldots,\gamma_{n};\alpha_{1},\beta_{1},\ldots,
\alpha_{g},\beta_{g}$ there is a bijective correspondence between the 
fiber of $\mathcal{H}_d^{(n)}(Y,b_0)\to (Y-b_0)^{(n)}-\Delta $ over 
$D$ and the set of Hurwitz systems 
$(t_1,\ldots,t_n;\lambda_1,\mu_1,\ldots,\lambda_g,\mu_g)$. Similarly 
there is a bijective correspondence between the fiber of 
$\mathcal{U}(b_{0})\to (Y-b_0)^{(n)}-\Delta $, where 
$\mathcal{U}(b_{0})\subset \mathcal{H}_d^{(n)}(Y)$, and the set of 
equivalence classes of Hurwitz systems 
$[t_1,\ldots,t_n;\lambda_1,\mu_1,\ldots,\lambda_g,\mu_g]$. In the next 
theorem we calculate the monodromy action of 
$\pi_1((Y-b_0)^{(n)}-\Delta ,D)$ on these fibers. According to 
Proposition~\ref{s2.19} it suffices to determine the braid moves 
which correspond to the generators $\sigma_{j},\rho_{ik},\tau_{ik}$. We 
denote the corresponding braid moves of the first type (cf. 
Definition~\ref{s2.17}) by $\sigma_{j}',\rho_{ik}',\tau_{ik}'$ and we 
denote the corresponding braid moves of the second type (inverse to those of 
the 
first type) by $\sigma_{j}'',\rho_{ik}'',\tau_{ik}''$.
\begin{thm}\label{s2.27}
Let $(t_1,\ldots,t_n;\lambda_1,\mu_1,\ldots,\lambda_g,\mu_g)$,
$\lambda_{g+k}=\mu_{k}$, be a Hurwitz system. Let 
$u_{k}=[\lambda_{1},\mu_{1}]\cdots [\lambda_{k},\mu_{k}]$ for 
$k=1,\ldots,g$ and let $u_{0}=1$. The following formulas hold for the 
braid moves 
$(t_1,\ldots,t_n;\lambda_1,\mu_1,\ldots,\lambda_g,\mu_g)\mapsto 
(t'_1,\ldots,t'_n;\lambda'_1,\mu'_1,\ldots,\lambda'_g,\mu'_g)$ of the 
first type.

\smallskip
a.\quad For $\sigma_{j}'$ where $1\leq j\leq n-1$
\begin{gather}
t'_{i}=t_{i}\quad \text{for $\forall i\neq j,\, j+1$},\quad
\lambda'_{\ell} = \lambda_{\ell},\quad \mu'_{\ell}=\mu_{\ell}\quad
\text{for $\forall \ell$}\notag\\
(t_{j},t_{j+1})\ \mapsto \ (t'_{j},t'_{j+1})\ =\ 
(t_{j}t_{j+1}t_{j}^{-1}\, ,\, t_{j}).\label{es3.27}
\end{gather}

\smallskip
b.\quad For $\rho_{ik}'$ where $1\leq i\leq n,\; 1\leq k\leq g$
\begin{gather}
t'_{j}=t_{j}\quad \text{for $\forall j\neq i$},\quad 
\lambda'_{\ell} = \lambda_{\ell}\quad \text{for $\forall \ell$},\quad 
\mu'_{\ell} = \mu_{\ell}\quad 
\text{for $\forall \ell \neq k$}\notag\\
(t_{i},\mu_{k})\ \mapsto \ (t'_{i},\mu'_{k})\ =\ 
(a_{1}^{-1}t_{i}a_{1},(b_{1}^{-1}t_{i}^{-1}b_{1})\mu_{k})\quad 
\text{where}\notag\\
a_{1}=(t_{1}\cdots 
t_{i-1})^{-1}u_{k-1}\lambda_{k}(u_{k}^{-1}u_{g})(t_{i+1}\cdots 
t_{n})^{-1},\quad b_{1}=(t_{1}\cdots t_{i-1})^{-1}u_{k-1}\lambda_{k}.\notag
\end{gather}

\smallskip
c.\quad For $\tau_{ik}'$ where $1\leq i\leq n,\; 1\leq k\leq g$
\begin{gather}
t'_{j}=t_{j}\quad \text{for $\forall j\neq i$},\quad 
\lambda'_{\ell} = \lambda_{\ell}\quad \text{for $\forall \ell \neq k$},
\quad 
\mu'_{\ell} = \mu_{\ell}\quad 
\text{for $\forall \ell$}\notag\\
(t_{i},\lambda_{k})\ \mapsto \ (t'_{i},\lambda'_{k})\ =\ 
(c_{1}^{-1}t_{i}c_{1},(d_{1}^{-1}t_{i}d_{1})\lambda_{k})\quad 
\text{where}\notag\\
c_{1}=t_{i+1}\cdots t_{n}(u_{k}^{-1}u_{g})^{-1}
\mu_{k}(u_{k-1})^{-1}t_{1}\cdots t_{i-1},\quad
d_{1}=t_{i+1}\cdots t_{n}(u_{k}^{-1}u_{g})^{-1}\mu_{k}. \notag
\end{gather}
The following formulas hold for the 
braid moves 
$(t_1,\ldots,t_n;\lambda_1,\mu_1,\ldots,\lambda_g,\mu_g)\mapsto 
\linebreak
(t''_1,\ldots,t''_n;\lambda''_1,\mu''_1,\ldots,\lambda''_g,\mu''_g)$ of the 
second type.

\smallskip
d.\quad For $\sigma_{j}''$ where $1\leq j\leq n-1$
\begin{gather}
t''_{i}=t_{i}\quad \text{for $\forall i\neq j,\, j+1$},\quad
\lambda''_{\ell} = \lambda_{\ell},\quad \mu''_{\ell}=\mu_{\ell}\quad
\text{for $\forall \ell$}\notag\\
(t_{j},t_{j+1})\ \mapsto \  (t''_{j},t''_{j+1})\ =\
(t_{j+1}\, ,\, t_{j+1}^{-1}t_{j}t_{j+1}).\label{es3.28a}
\end{gather}

\smallskip
e.\quad For $\rho_{ik}''$ where $1\leq i\leq n,\; 1\leq k\leq g$
\begin{gather}
t''_{j}=t_{j}\quad \text{for $\forall j\neq i$},\quad 
\lambda''_{\ell} = \lambda_{\ell}\quad \text{for $\forall \ell$},\quad 
\mu''_{\ell} = \mu_{\ell}\quad 
\text{for $\forall \ell \neq k$}\notag\\
(t_{i},\mu_{k})\ \mapsto \ (t''_{i},\mu''_{k})\ =\ 
(a_{2}^{-1}t_{i}a_{2},(b_{2}^{-1}t_{i}b_{2})\mu_{k})\quad 
\text{where}\notag\\
a_{2}= t_{i+1}\cdots t_{n}(u_{k}^{-1}u_{g})^{-1}\lambda_{k}^{-1}
(u_{k-1})^{-1}t_{1}\cdots t_{i-1},\quad 
b_2 = t_{i+1}\cdots t_{n}(u_{k}^{-1}u_{g})^{-1}.\notag
\end{gather}

\smallskip
f.\quad For $\tau_{ik}''$ where $1\leq i\leq n,\; 1\leq k\leq g$
\begin{gather}
t''_{j}=t_{j}\quad \text{for $\forall j\neq i$},\quad 
\lambda''_{\ell} = \lambda_{\ell}\quad \text{for $\forall \ell \neq k$},
\quad 
\mu''_{\ell} = \mu_{\ell}\quad 
\text{for $\forall \ell$}\notag\\
(t_{i},\lambda_{k})\ \mapsto \ (t''_{i},\lambda''_{k})\ =\ 
(c_{2}^{-1}t_{i}c_{2},(d_{2}^{-1}
t_{i}^{-1}d_{2})\lambda_{k})\quad 
\text{where}\notag\\
c_{2}=(t_{1}\cdots t_{i-1})^{-1}u_{k-1}\mu_{k}^{-1}
(u_{k}^{-1}u_{g})(t_{i+1}\cdots t_{n})^{-1},\quad 
d_{2}=(t_{1}\cdots t_{i-1})^{-1}u_{k-1}. \notag
\end{gather}
\end{thm}
\begin{proof}
Formulas (d), (e) and (f) are obtained from formulas (a), 
(b) and (c) of Theorem~\ref{s2.20} respectively  applying the homomorphism 
$m:\pi_1(Y-D,b_{0})\to S_{d}$ and equalities Eq.\eqref{es2.16c} and 
Eq.\eqref{es2.16b}. By Lemma~\ref{s2.16} the braid moves of the first 
type are inverse to the braid moves of the second type. We  may thus
obtain formulas (a), (b) and (c) from formulas (d), (e) 
and (f) of Theorem~\ref{s2.20} respectively  applying the homomorphism $m$. 
\end{proof}
\begin{cor}\label{s2.30a}
Using the notation of Theorem~\ref{s2.27} the following formulas hold.
\renewcommand{\theenumi}{\roman{enumi}}
\begin{enumerate}
\item For $\rho_{ik}'$\; :
\begin{alignat*}{2}
&\rho'_{nk}:t_{n}\mapsto t'_{n} = a_{1}^{-1}t_{n}a_{1}, 
&\quad
&\text{where
$a_{1}=(u_{k}^{-1}u_{g})^{-1}[\mu_{k},\lambda_{k}]\lambda_{k}
(u_{k}^{-1}u_{g})$}\\
&\rho_{1k}':\mu_{k}\mapsto 
\mu_{k}'=(b_{1}^{-1}t_{1}^{-1}b_{1})\mu_{k},
&
&\text{where
$b_{1}=u_{k-1}\lambda_{k}$}\\
&\rho_{nk}':\mu_{k}\mapsto 
\mu'_{k}=(b_{1}^{-1}t_{n}^{-1}b_{1})\mu_{k},
&
&\text{where $b_{1}=
(u_{k}^{-1}u_{g})^{-1}[\mu_{k},\lambda_{k}]\lambda_{k}$}.
\end{alignat*}
\item For $\tau'_{ik}$\; :
\begin{align*}
&\tau'_{1k}:t_{1}\mapsto t'_{1} = c_{1}^{-1}t_{1}c_{1}, 
&\quad 
&\text{where 
$c_{1}=
u_{k-1}[\lambda_{k},\mu_{k}]\mu_{k}u_{k-1}^{-1}$\hspace{0.3cm}}\\
&\tau'_{1k}:\lambda_{k}\mapsto 
\lambda_{k}'=(d_{1}^{-1}t_{1}d_{1})\lambda_{k},
& 
&\text{where $d_{1}=u_{k-1}[\lambda_{k},\mu_{k}]\mu_{k}$}\\
&\tau_{nk}':\lambda_{k}\mapsto \lambda_{k}'=
(d_{1}^{-1}t_{n}d_{1})\lambda_{k},
& 
&\text{where $d_{1}=(u_{k}^{-1}u_{g})^{-1}\mu_{k}$}.
\end{align*}
\item For $\rho''_{ik}$\; : 
\begin{align*}
&\rho''_{1k}:t_{1}\mapsto t''_{1} = a_{2}^{-1}t_{1}a_{2}, 
&\quad
&\text{where
$a_{2}=u_{k-1}[\lambda_{k},\mu_{k}]\lambda_{k}^{-1}
u_{k-1}^{-1}$}\\
&\rho_{1k}'':\mu_{k}\mapsto 
\mu_{k}''=(u_{k}^{-1}t_{1}u_{k})\mu_{k},
&
&\\
&\rho_{nk}'':\mu_{k}\mapsto 
\mu''_{k}=(b_{1}^{-1}t_{n}b_{1})\mu_{k},
&
&\text{where $b_{1}=
(u_{k}^{-1}u_{g})^{-1}$}.
\end{align*}
\item For $\tau''_{ik}$\; :
\begin{alignat*}{2}
&\tau''_{nk}:t_{n}\mapsto t''_{n} = c_{2}^{-1}t_{n}c_{2}, 
&\quad 
&\text{where 
$c_{2}=
(u_{k}^{-1}u_{g})^{-1}[\mu_{k},\lambda_{k}]\mu_{k}^{-1}
(u_{k}^{-1}u_{g})$}\\
&\tau''_{1k}:\lambda_{k}\mapsto 
\lambda_{k}''=(u_{k-1}^{-1}t_{1}^{-1}u_{k-1})\lambda_{k},
& 
&\\
&\tau_{nk}'':\lambda_{k}\mapsto \lambda_{k}''=
(d_{2}^{-1}t_{n}^{-1}d_{2})\lambda_{k},
& 
&\text{where $d_{2}=(u_{k}^{-1}u_{g})^{-1}[\mu_{k},\lambda_k]$}.
\end{alignat*}
\end{enumerate}
In particular
\begin{equation}\label{es2.30b}
\rho_{ng}''\ :\ \mu_{g}\ \mapsto \ t_n\mu_g,\quad \quad \tau''_{11}\ :\ 
\lambda_{1}\ \mapsto \ t_1^{-1}\lambda_{1}.
\end{equation}
\end{cor}
\begin{proof}
Let us prove the first formula of Part~(ii). The other formulas can be 
either proved similarly or are restatements of particular cases of 
Theorem~\ref{s2.27}. We have $\tau_{1k}':t_{1}\mapsto 
t'_1=c_1^{-1}t_1c_1$ where $c_1=t_{2}\cdots 
t_{n}(u_k^{-1}u_g)^{-1}\mu_{k}u_{k-1}^{-1}$. Since $t_{1}\cdots 
t_{n}=u_g$ it holds $(t_{2}\cdots t_{n})^{-1}t_1(t_{2}\cdots 
t_{n})=u_g^{-1}t_1u_g$. Hence
\begin{equation*}
\begin{split}
t'_{1}\ 
&= \ (u_{k-1}\mu_{k}^{-1}u_k^{-1})t_1(\cdots )^{-1}\ =\ 
(u_{k-1}\mu_{k}^{-1}[\lambda_{k},\mu_{k}]^{-1}
u_{k-1}^{-1})t_1(\cdots)^{-1}\\
&= \ (u_{k-1}[\lambda_{k},\mu_{k}]\mu_{k}u_{k-1}^{-1})^{^{-1}}\, t_1\,
(u_{k-1}[\lambda_{k},\mu_{k}]\mu_{k}u_{k-1}^{-1})
\end{split}
\end{equation*}
\end{proof}
Summing up the discussion made so far in this section we obtain the 
following result.
\begin{thm}\label{s2.30c}
Fix $d\geq 2,\; n>0,\; g\geq 0$. Let us consider the set of all 
Hurwitz systems 
$(t_1,\ldots,t_n;\lambda_1,\mu_1,\ldots,\lambda_g,\mu_g)$ (cf. 
Definition~\ref{s2.14}). Let $F$ be the free group generated by the 
symbols $\sigma_{j},\rho_{ik},\tau_{ik}$ where $1\leq j\leq n-1,\; 
1\leq i\leq n,\; 1\leq k\leq g$. Let us consider the right action of 
$F$ on the set of Hurwitz systems defined by the formulas for 
$\sigma_{j}',\rho_{ik}',\tau_{ik}'$ of Theorem~\ref{s2.27}~(a)--(c) or 
alternatively let us consider the left action of $F$ on the set of 
Hurwitz systems defined by the formulas for 
$\sigma_{j}'',\rho_{ik}'',\tau_{ik}''$ of Theorem~\ref{s2.27}~(d)--(f). 
Then the connected components of $\mathcal{H}_d^{(n)}(Y,b_0)$ (cf. 
\ref{s2.11}) correspond bijectively to the orbits of either action of 
$F$ and the connected components of $\mathcal{H}_d^{(n)}(Y)$ 
correspond bijectively to the orbits of either associated action of 
$F$ on the set of equivalence classes of Hurwitz systems modulo inner 
automorphisms of $S_d$. 
\end{thm}
\begin{proof}
$\mathcal{H}_d^{(n)}(Y,b_0)\to (Y-b_0)^{(n)}-\Delta $ is a topological 
covering map where
a Hausdorff topology on $\mathcal{H}_d^{(n)}(Y,b_0)$ is defined as in 
\cite{Fu} pp.545,546. According to the definition of product of arcs 
Eq.\eqref{es2.12} the monodromy action of the fundamental group 
$\pi_1((Y-b_0)^{(n)}-\Delta ,D)$ on the fiber over $D$
is a right action. The identification of this fiber with the set of Hurwitz 
systems, Proposition~\ref{s2.19} and the calculation of the braid moves of 
the first type $\sigma_{j}',\rho_{ik}',\tau_{ik}'$ in 
Theorem~\ref{s2.27} yield the statement about the orbits of the right 
action of $F$.

Consider the associated left action $gx=xg^{-1}$. The orbits are the 
same and according to Lemma~\ref{s2.16} the braid moves of the second 
type $\sigma_{j}'',\rho_{ik}'',\tau_{ik}''$ are inverse to 
 those of the first type 
$\sigma_{j}',\rho_{ik}',\tau_{ik}'$ respectively. This shows the statement about 
the orbits of the left action of $F$. 

Let $b_{0}\in Y$. Since $\mathcal{U}(b_{0})\subset 
\mathcal{H}_d^{(n)}(Y)$ (cf. \ref{s2.11}) is a Zariski dense open 
subset, the connected components of $\mathcal{U}(b_{0})$ correspond 
bijectively to those of $\mathcal{H}_d^{(n)}(Y)$. This shows the last 
claim of the theorem.
\end{proof}

\section{The main lemma}\label{s3}
In the previous section we considered Hurwitz systems 
$(t_1,\ldots,t_n;\lambda_1,\mu_1,\ldots,\lambda_g,\mu_g)$ where 
$t_i,\lambda_{k},\mu_{k}\in S_d$ in connection with the problem of 
determining the connected components of the Hurwitz spaces. In this 
section we replace $S_d$ by an arbitrary (possibly infinite) group 
$G$. 
\begin{dfn}\label{s3.31a}
Let $G$ be an arbitrary group. We call 
$(t_1,\ldots,t_n;\lambda_1,\mu_1,\ldots,\lambda_g,\mu_g)$ a Hurwitz 
system with values in $G$ if $t_i,\lambda_{k},\mu_{k}\in G$, $t_i\neq 
1$ for $\forall i$ and $t_{1}\cdots t_{n}=[\lambda_{1},\mu_{1}]\cdots 
[\lambda_{g},\mu_{g}]$.
\end{dfn}
The formulas of Theorem~\ref{s2.27} make sense for an arbitrary group 
$G$. That $\sigma_{j}''=(\sigma_{j}')^{-1},\: 
\rho_{ik}''=(\rho_{ik}')^{-1}$ and $\tau_{ik}''=(\tau_{ik}')^{-1}$ is 
evident from Theorem~\ref{s2.20}.
\begin{dfn}\label{s3.31ab}
We call two Hurwitz systems with values in $G$ braid-equivalent if one 
is obtained from the other by a finite sequence of braid moves 
$\sigma_{j}',\rho_{ik}',\tau_{ik}',\sigma_{j}'',\rho_{ik}'',
\tau_{ik}''$ where $1\leq j\leq n-1,\: 1\leq i\leq n,\: 1\leq k\leq 
g$. We denote the braid equivalence by $\sim$. 
\end{dfn}
\begin{mainlem}\label{s3.32}
Let $G$ be an arbitrary group and let 
$(t_1,\ldots,t_n;\lambda_1,\mu_1,\ldots,\lambda_g,\mu_g)$ be a Hurwitz 
system with values in $G$. Suppose that $t_it_{i+1}=1$. Let $H$ be the 
subgroup of $G$ generated by 
$\{t_1,\ldots,t_{i-1}, 
t_{i+2},\ldots, t_n;\lambda_1,\mu_1,\ldots,\lambda_g,\mu_g\}$. Then for 
every $h\in H$ the given Hurwitz system is braid-equivalent to 
\(
(t_1,\ldots,t_{i-1},t_{i}^{h},t_{i+1}^{h},
t_{i+2},\ldots, t_n;\lambda_1,\mu_1,\ldots,\lambda_g,\mu_g)
\)
\end{mainlem}
\begin{proof}
Let us fix 
\(
t_1,\ldots,t_{i-1}, 
t_{i+2},\ldots, t_n;\lambda_1,\mu_1,\ldots,\lambda_g,\mu_g
\). Let $H_1\subseteq H$ be the subset consisting of elements $h$
such that the statement of the lemma holds for an arbitrary pair 
$(t_i,t_{i+1})=(\tau,\tau^{-1})$. 
\begin{stepa}
\emph{We claim $H_1$ is a subgroup of $H$}. Let $h\in H_1$ and let 
$t=ht_ih^{-1}$. Then $t_i=t^{h}$, so by assumption 
$(\ldots,t^{h},(t^{-1})^{h},\ldots)$ can be obtained from 
$(\ldots,t,t^{-1},\ldots)$ by a sequence of braid moves. Then one can 
obtain $(\ldots,t,t^{-1},\ldots)$ from $(\ldots,t_i,
t_i^{-1},\ldots)$ 
by the inverse sequence of braid moves. Thus $h^{-1}\in H_1$. If 
$h_1,h_2\in H_1$, then 
\(
(\ldots,t_i,t_i^{-1},\ldots)\sim 
(\ldots,t_i^{h_1},(t_i^{-1})^{h_1},\ldots)\sim
(\ldots,t_i^{h_1h_2},(t_i^{-1})^{h_1h_2},\ldots)
\),
so $h_1h_2\in H_1$.
\end{stepa}
\begin{stepa}
\emph{For every $\ell\neq i,i+1$ the element $t_{\ell}$ belongs to $H_1$}. 
Applying a sequence of braid moves $\sigma_{j}',\sigma_{j}''$ we can 
move the adjacent pair $(t_i,t_{i+1})$ wherever we want among the 
first $n$ elements of the Hurwitz system without changing the other 
elements. So move $(t_i,t_{i+1})$ to the left side of $t_{\ell}$. Then 
we have
\(
(t_{i},t_{i+1},t_{\ell})\sim 
(t_{i},t_{\ell},t_{\ell}^{-1}t_{i+1}t_{\ell})\sim 
(t_{\ell},t_{\ell}^{-1}t_{i}t_{\ell},t_{\ell}^{-1}t_{i+1}t_{\ell})
\).
We then move the pair $(t_i^{h},t_{i+1}^{h})$ with $h=t_{\ell}$ back 
to the initial position.
\end{stepa}
\begin{stepa}\label{stepa3}
\emph{For every $k=1,\ldots,g$ the element 
$h=u_{k-1}\lambda_{k}u_{k}^{-1}$ belongs to $H_1$}. First suppose that $i=1$. In this case
$t_{2}=t_{1}^{-1}$. Let us perform a braid move $\rho_{1k}'$. One 
obtains 
\(
(t_1,t_1^{-1},\ldots,\lambda_k,\mu_k,\ldots) \sim 
(t'_1,t_1^{-1},\ldots,\lambda_k,\mu'_k,\ldots)
\)
where $t_1'=a_1^{-1}t_1a_1$ with 
$a_1=\linebreak
u_{k-1}\lambda_{k}(u_k^{-1}u_g)(t_{2}\cdots t_{n})^{-1}$ and 
$\mu_{k}'=(b_1^{-1}t_1^{-1}b_1)\mu_{k}$ with $b_1=u_{k-1}\lambda_{k}$. 
We have $t_{1}\cdots t_{n}=u_g$, so 
$a_1=u_{k-1}\lambda_{k}u_k^{-1}t_1$. Let us move $t_2=t_1^{-1}$ to the 
first place using $\sigma_{1}''$. One obtains 
$(t_1^{-1},t_1t'_1t_1^{-1},\ldots,\lambda_k,\mu'_k,\ldots)$. Here 
$t_1t_1't_1^{-1}=t_1^{h}$ where $h=u_{k-1}\lambda_{k}u_{k}^{-1}$. Let 
us perform again $\rho_{1k}'$. The elements $u_{k-1},\lambda_{k}$ have 
not changed with respect to the original Hurwitz system, so 
$\mu_{k}'\mapsto 
[(u_{k-1}\lambda_{k})^{-1}(t_1^{-1})^{-1}(u_{k-1}\lambda_{k})]\mu_{k}' 
= \mu_{k}$. Moving the second element $t_1^{h}$ to the first place 
by $\sigma_{1}''$ we obtain  
$(t_1^{h},\tilde{t}_{2},t_3,\ldots,\lambda_{k},\mu_{k},\ldots)$. Since 
$t_{1}t_{2}t_{3}\cdots t_{n}=u_{g}=t_{1}^{h}\tilde{t}_2t_3\cdots 
t_{n}$ we must have $\tilde{t}_2=(t_1^{-1})^{h}$. This proves 
$(t_1,t_1^{-1},t_3,\ldots)\sim (t_1^{h},(t_1^{-1})^{h},t_3,\ldots)$ 
with $h=u_{k-1}\lambda_{k}u_{k}^{-1}$. One extends this braid 
equivalence to every adjacent pair $(t_i,t_{i+1})$ with $t_it_{i+1}=1$ 
by moving first the pair to the front, applying the braid equivalence 
we have just proved and moving the obtained pair back to the original 
position. 
\end{stepa}
\begin{stepa}
\emph{For every $k=1,\ldots,g$ the element 
$h=u_{k-1}\mu_{k}^{-1}u_k^{-1}$ belongs to $H_1$}. The proof is the 
same as that of Step~\ref{stepa3}. One uses the braid move 
$\tau_{1k}''$ instead of $\rho_{1k}'$. 
\end{stepa}
\begin{stepa}
By the preceding steps it remains to verify that 
$\lambda_{k},\mu_{k}$ for $\forall k$ belong to the subgroup 
$H_2\subseteq H_1$ generated by $t_j,
u_{k-1}\lambda_{k}u_{k}^{-1},
u_{k-1}\mu_{k}^{-1}u_k^{-1}$ where $j=1,\ldots,i-1,i+2,\ldots,n$ and 
$k=1,\ldots,g$. We prove this by induction on $k$. We have 
$\lambda_{1}u_1^{-1}\in H_2,\; 
u_1\mu_{1}=(\mu_{1}^{-1}u_1^{-1})^{-1}\in H_2$, so 
$\lambda_{1}\mu_{1}\in H_2$. Furthermore $H_2\ni 
u_1\mu_{1}=[\lambda_{1},\mu_{1}]\mu_{1}=\lambda_{1}\mu_{1}
\lambda_{1}^{-1}$. So $\lambda_{1}^{^{\pm 1}}\in H_2$ and 
$\mu_{1}^{^{\pm 1}}\in H_2$. Let $k\geq 2$. Suppose, by inductive
assumption, that $\lambda_1,\mu_1,\ldots,\lambda_{k-1},\mu_{k-1}$ belong to 
$H_{2}$. Then $u_{k-1}=\prod_{\ell 
=1}^{k-1}[\lambda_{\ell},\mu_{\ell}]\in H_2$. We have 
$u_{k-1}\lambda_{k}u_{k}^{-1}\in H_2$, 
$u_k\mu_{k}u_{k-1}^{-1} = (u_{k-1}\mu_{k}^{-1}u_k^{-1})^{-1}\in H_2$, 
so $u_{k-1}\lambda_{k}\mu_{k}u_{k-1}^{-1}\in H_2$ and therefore 
$\lambda_{k}\mu_{k}\in H_2$. Furthermore $H_2\ni 
u_k\mu_{k}u_{k-1}^{-1}=u_{k-1}[\lambda_{k},\mu_{k}]\mu_{k}u_{k-1}^{-1}
$. Therefore $\lambda_{k}\mu_{k}\lambda_{k}^{-1}\in H_2$, so 
$\lambda_{k}^{^{\pm 1}}\in H_2$ and 
$\mu_{k}^{^{\pm 1}}\in H_2$. 
\end{stepa}
The lemma is proved.
\end{proof}

\section{The case $n\geq 2d-2$}\label{s4}
So far we have not made any restrictions on the type of the covering 
$\pi :X\to Y$. We shall be further occupied mainly with simply branched 
coverings with connected $X$. These coverings correspond to 
Hurwitz systems 
$(t_1,\ldots,t_n;\lambda_1,\mu_1,\ldots,\lambda_g,\mu_g)$ with local 
monodromies $t_i=m(\gamma_{i})$, $i=1,\ldots,n$ equal to 
transpositions and with transitive monodromy group $G=\langle 
t_1,\ldots,t_n,\lambda_1,\mu_1,\ldots,\lambda_g,\mu_g
\rangle$. We call such coverings \emph{simple coverings}. By Hurwitz' 
formula $n\equiv 0(mod\; 2)$ for a simple covering. 
\begin{dfn}\label{s4.37}
Let $Y$ be a smooth projective curve. We denote by 
$\mathcal{H}_{d,n}(Y,b_0)$ and $\mathcal{H}_{d,n}(Y)$ the Hurwitz 
spaces which parametrize equivalence classes $[X\to Y,\: \phi]$ and  
$[X\to Y]$ respectively of simple coverings of degree $d$ branched in $n$ points 
(cf. \ref{s2.11}). 
\end{dfn}
As a first application of Theorem~\ref{s2.30c} we give a proof, using 
braid moves, of the following well-known fact.
\begin{pro}\label{s4.37aa}
Let $n>0,\; n\equiv 0(mod\: 2)$. The Hurwitz spaces 
$\mathcal{H}_{2,n}(Y,b_0)$ and $\mathcal{H}_{2,n}(Y)$ which 
parametrize double coverings of $Y$ branched in $n$ points are both 
irreducible.
\end{pro}
\begin{proof}
Since $\mathcal{H}_{2,n}(Y,b_0)$ and $\mathcal{H}_{2,n}(Y)$ are smooth 
it suffices to prove they are connected. Let us first consider 
$\mathcal{H}_{2,n}(Y,b_0)$. By Theorem~\ref{s2.30c} it suffices to 
prove that every Hurwitz system 
$(t_1,\ldots,t_n;\lambda_1,\mu_1,\ldots,\lambda_g,\mu_g)$
with $t_i,\lambda_{k},\mu_{k}\in S_2$ is braid-equivalent to 
$((12),\ldots,(12);1,1,\ldots,\linebreak
1,1)$. According to Theorem~\ref{s2.27} 
one has $\rho_{ik}':(t_i,\mu_{k})\mapsto (t_i,t_i\mu_{k})$, 
$\tau_{ik}':(t_i,\lambda_{k})\mapsto (t_i,t_i\lambda_{k})$ and the 
same formulas hold for $\rho_{ik}''$ and $\tau_{ik}''$. Hence whenever 
some $\lambda_{k}\neq 1$ or $\mu_{k}\neq 1$ we can perform a braid move 
$\tau_{1k}'$ or $\rho_{1k}'$ respectively or both in order to obtain a new Hurwitz 
system with $\lambda_{k}'=\mu_{k}'=1$. The connectedness of 
$\mathcal{H}_{2,n}(Y,b_0)$ implies the connectedness of the Zariski dense 
open subset $\mathcal{U}(b_0)\subset \mathcal{H}_{2,n}(Y)$. Therefore 
$\mathcal{H}_{2,n}(Y)$ is connected as well.
\end{proof}
Given an ordered $n$-tuple of permutations $\underline{t}=
(t_1,\ldots,t_n)$
 whose product $t_{1}\cdots t_{n}=s\in S_d$, performing an elementary 
move $\sigma_{j}'$ or $\sigma_{j}''$ one obtains a new ordered 
$n$-tuple $(t_1',\ldots,t_n')$ which has the same product $t_{1}'\cdots 
t_{n}'=s$. We shall also call ordered $n$-tuples sequences.
\begin{dfn}\label{s4.37ab}
Two ordered $n$-tuples (or sequences) of permutations 
$\underline{t}=(t_1,\ldots,t_n)$ and 
$\underline{t}'=(t_1',\ldots,t_n')$ are called braid-equivalent if 
$\underline{t}'$ is obtained from $\underline{t}$ by a finite sequence 
of braid moves of type $\sigma_{j}'$ or $\sigma_{j}''$. We write 
$\underline{t}'\sim \underline{t}$. 
\end{dfn}
Notice the difference between braid equivalence of sequences and that 
of Hurwitz systems (cf. Definition~\ref{s3.31ab}). The latter includes 
also braid moves of the types 
$\rho_{ik}',\rho_{ik}'',\tau_{ik}',\tau_{ik}''$.
\begin{block}\label{s4.37b}
We recall some results proved in \cite{Mo}. Given a permutation $s\in 
S_d$ one considers an ordered $n$-tuple of transpositions 
$\underline{t}=(t_1,\ldots,t_n)$ such that $t_{1}\cdots t_{n}=s$. 
Let $\Sigma_1,\Sigma_2,\ldots,\Sigma_m$ be the domains of transitivity 
of the group $G=\langle t_1,\ldots,t_n\rangle$. Then 
$G=S(\Sigma_1)\times\cdots \times S(\Sigma_m)$. For reader's convenience we 
give a proof of the following important lemma due to Mochizuki (cf. 
\cite{Mo} Lemma~2.4).
\end{block}
\begin{lem}\label{s4.37a}
Let $\underline{t}=(t_1,\ldots,t_n),\; t_{1}\cdots t_{n}=s$ be as 
above. Let $a,b\in \Sigma_i,\; a<b$. Then $\underline{t}\sim 
\underline{t}'=((ab),t'_2,\ldots)$ and 
$\underline{t}\sim \underline{t}''=(\ldots,t''_{n-1},(ab))$.
\end{lem}
\begin{proof}
If $(t_i,t_{i+1})$ is a pair such that $t_i\in S(\Sigma_{k}),\; 
t_{i+1}\in S(\Sigma_{\ell})$ with $k\neq \ell$, then $(t_i,t_{i+1})\sim 
(t_{i+1},t_{i+1}^{-1}t_it_{i+1})=(t_{i+1},t_i)$. This shows that 
performing a 
sequence of braid moves (of type $\sigma_{j}',\sigma_{j}''$) one may 
replace $\underline{t}$ by a concatenation of sequences $T_1T_2\ldots 
T_m$ where $T_i$ is formed by all transpositions of $\underline{t}$ 
which belong to $S(\Sigma_j)$ ordered in the way they appear in 
$\underline{t}$ (if $\#\, \Sigma_{j}=1$ we let $T_j=\emptyset$). 
Furthermore we may move any $T_i$ to the first place or to the last place. 
This shows that it 
suffices to prove the lemma in the particular case when $\langle 
t_{1},\ldots ,t_{n} \rangle$ is a transitive group. Let $b=a^{\tau}$ 
where $\tau=t_{j_1}t_{j_2}\cdots t_{j_{\ell}}$. One may vary the 
sequence within the set of sequences braid-equivalent to 
$\underline{t}$. Given a sequence one may vary $\tau$ so that 
$b=a^{\tau}$ and finally given a sequence and a $\tau$ one may vary 
the factorization of $\tau$. Let $\underline{t}'$, $\tau'$ with the property 
$a^{\tau'}=b$ and a factorization 
$\tau'=t'_{i_1}\cdots t'_{i_r}$ be chosen so that $r$ is minimal 
possible. If $r=1$ one has that $(ab)$ is one of the transpositions of 
$\underline{t}'$. Moving $(ab)$ to the front by subsequent 
elementary moves $\sigma_{j}''$ one obtains a braid-equivalent 
sequence of the type $((ab),\ldots)$. Moving $(ab)$ to the end by 
subsequent 
elementary moves $\sigma_{j}'$ one obtains $(\ldots,(ab))$ as required 
in the lemma. Suppose that $r\geq 2$. Let $x_1=a,\, x_2=a^{t'_{i_1}},\ldots, 
x_{k+1}=a^{t'_{i_1}\cdots t'_{i_k}},\ldots, x_{r+1}=b$. The 
minimality of $r$ implies that $x_i\neq x_j$ for $i\neq j$, so 
$t'_{i_k}=(x_k\, x_{k+1})$ for $k=1,\ldots,r$ are $r$ different 
transpositions of the sequence $\underline{t}'$ and every $x_i$ for 
$2\leq i\leq r$ enters in at most two transpositions of the set  
$\{t'_{i_k}\}$. Applying several 
elementary moves one 
places $(ax_2)$ adjacent to $(x_2 x_3)$. One obtains a sequence either 
of the type $\ldots,(a\, x_2),(x_2x_3),\ldots$ or of the type 
$\ldots,(x_2x_3),(a\, x_2),\ldots$ in which none of $t'_{i_k},\, 
k=1,\ldots,r$ has been changed. In the first case one has 
$((a\, x_2),(x_2x_3))\sim ((x_2x_3),(a\, x_3))$ and in the second case 
$((x_2x_3),(a\, x_2))\sim ((a\, x_3),(x_2x_3))$. In both cases one 
obtains a braid-equivalent sequence for which the sequence 
$a=x_1,x_2,x_3,\ldots,x_{r+1}=b$ is replaced by 
$a=x_1,x_3,\ldots,x_{r+1}=b$. This contradicts the minimality of $r$. 
\end{proof}
\begin{block}\label{s4.40a}
Let $(t_1,\ldots,t_n)$ be an $n$-tuple of transpositions and let 
$t_{1}\cdots t_{n}=s$. Let $\Delta_{1},\ldots,\Delta_{q}$ be the 
domains of transitivity of the cyclic group $\langle s\rangle$ where 
$\#\, \Delta_{i}=e_i\geq 1$. Following \cite{Mo} pp.369,370 if 
$s=s_1\cdots s_q$ is the corresponding product of independent cycles 
we may write $s_i=(1_i2_i\ldots (e_i)_i)$ if $e_i>1$ and $s_i=(1_i)=1$ if 
$e_i=1$. Such a representation is uniquely determined by $s$ if we 
assume that that $1_i$ is the minimal number among $\Delta_{i}$ for each 
$i$ and if we order $\Delta_{i}$ in such a way that $1_1<1_2<\cdots <
1_q$. Let $e_i>1$ and let $\mathcal{Z}_{i}$ be the sequence 
$((1_i2_i),(1_i3_i),\ldots,(1_i(e_i)_i))$. If $e_i=1$ one lets 
$\mathcal{Z}_{i}=\emptyset$. Let $\mathcal{Z}$ be the concatenation of 
sequences $\mathcal{Z}=\mathcal{Z}_{1}\mathcal{Z}_{2}\ldots 
\mathcal{Z}_{q}$. The sequence $\mathcal{Z}$ consists of 
$N=\sum_{i=1}^{q}(e_i-1)$ transpositions whose product equals $s$. The 
following theorem is proved in \cite{Mo} pp.369,370. It can also be 
deduced from the earlier paper of Kluitmann \cite{Kl}.
\end{block}
\begin{thm}\label{s4.40}
Given a sequence of transpositions $\underline{t}=(t_1,\ldots,t_n)$ 
satisfying $t_{1}\cdots t_{n}=s$ and such that $\langle t_{1},\ldots 
,t_{n} \rangle$ is a transitive group (therefore = $S_d$), there is a 
braid-equivalent sequence $(t'_1,\ldots,t'_n)$ which has the following 
form
\[
(t'_1,\ldots,t'_n)\ =\ (\mathcal{Z},t'_{N+1},\ldots,t'_n)
\]
where $n-N\equiv 0(mod\, 2)$ and 
\begin{enumerate}
\item if $q=1$, then $t'_i=(1_12_1)$ for $\forall i\geq N+1$
\item if $q>1$ then
\[
(t'_{_{N+1}},\ldots,t'_n)\ =\ 
((1_11_2),(1_11_2),(1_11_3),(1_11_3),\ldots ,(1_11_q),(1_11_q))
\]
where each $(1_11_i)$ appears two times if $2\leq i\leq q-1$ and 
$(1_11_q)$ appears an even number of times .
\end{enumerate}
\end{thm}
\begin{cor}\label{s4.41}
Let $\underline{t}=(t_1,\ldots,t_n)$ be a sequence of transpositions 
in $S_d$ which generate a transitive group. Let $s=t_{1}\cdots t_{n}$ 
and let $\Delta_{1},\ldots,\Delta_{q}$ be the domains of transitivity 
of the cyclic group $\langle s\rangle$. Suppose that $q\geq 2$. Let $a\in 
\Delta_{i},\; b\in \Delta_{j}$ where $i\neq j$. Then $\underline{t}$ 
is braid-equivalent to a sequence $(\ldots,(ab),(ab))$. 
\end{cor}
\begin{proof}
Let us consider the group $S_q$ which permutes the indices of $\Delta_{i},\, 
i=1,\ldots,q$. 
Let $\tau_1=(ij),\tau_{2},\ldots,\tau_{q-1}$ be transpositions in 
$S_q$ which generate it. We let $a_i=a,\; a_j=b$ and choose 
arbitrarily $a_k\in \Delta_{k}$ for each $k\neq i,j$. If 
$\tau_{\ell}=(\alpha,\beta)$ we consider the pair of transpositions 
$(a_{\alpha}a_{\beta}),(a_{\alpha}a_{\beta})$. We then consider a 
sequence of $n$ transpositions in $S_d$ which is a concatenation of the 
sequence
$\mathcal{Z}$, of the pairs associated with 
$\tau_{1},\tau_{2},\ldots,\tau_{q-1}$ and of the sequence $(1\, 2),\ldots,(1\, 2)$ with $(12)$ repeated 
an even number of times. We obtain a sequence 
$\mathcal{Z},(ab),(ab),\ldots$ with product $s$ which generates a 
transitive group. According to Theorem~\ref{s4.40} this sequence is 
braid-equivalent to $\underline{t}$. Moving the pair $(ab),(ab)$ to the end we 
obtain 
the required sequence.
\end{proof}
\begin{block}\label{s4.41a}
Given a permutation $s\in S_d$ whose domains of transitivity are 
$\Delta_{1},\ldots,\Delta_{q}$ with $\# \, \Delta_{i}=e_i$ we let 
$|s|=\sum_{i=1}^{q}(e_i-1)$. Let us consider a sequence of 
transpositions $\underline{t}=(t_1,\ldots,t_n),\; t_i\in S_d$. Let 
$\underline{\Sigma}=\{\Sigma_1,\ldots,\Sigma_m\}$ be the domains of 
transitivity of the group $\langle t_{1},\ldots ,t_{n} \rangle$. We 
let $|\underline{\Sigma}|=\sum_{j=1}^{m}(\# \,\Sigma_{j}-1)$.
\end{block}
\begin{lem}\label{s4.42}
Let $\underline{t}=(t_1,\ldots,t_n)$ be a sequence of transpositions 
and let $t_{1}\cdots t_{n}=s$. Then $n+|s|\equiv 0(mod\, 2)$ and 
$n+|s|\geq 2|\underline{\Sigma}|$. If 
$n+|s|> 2|\underline{\Sigma}|$ then there is a sequence 
$\underline{t}'$ braid-equivalent to $\underline{t}$ such that 
$t'_{n-1}=t'_{n}$ and $t'_1,\ldots,t'_{n-2}$ generate the same group 
as $t_1,\ldots,t_n$.
\end{lem}
\begin{proof}
It is clear that every $\Delta_{i}$ belongs to some $\Sigma_j$. 
Furthermore one may replace $\underline{t}$ by a braid-equivalent 
sequence which is a concatenation $T_1T_2\ldots T_m$ of sequences 
$T_j$ whose elements belong to $S(\Sigma_j)$ (cf. the proof of 
Lemma~\ref{s4.37a}; if $\# \, \Sigma_j=1$ one lets $T_j=\emptyset$). It 
thus suffices to prove the lemma in the case when $\langle 
t_{1},\ldots ,t_{n} \rangle$ is transitive, which we further assume, so 
$\underline{\Sigma}=\{\Sigma_1\},\; |\underline{\Sigma}|=d-1$. Let us 
consider the braid-equivalent sequence $\underline{t}'$ of 
Theorem~\ref{s4.40}. The subsequence $\mathcal{Z}$ has $N=|s|$ 
elements. The remaining $n-|s|$ elements appear in pairs and their 
number is at least $\geq 2(q-1)$ since $\langle t'_{1},\ldots ,t'_{n} 
\rangle$ is a transitive group. If $n-|s|>2(q-1)$ the pair 
$(1_11_q),(1_11_q)$ is repeated at least twice, so canceling it from 
the sequence does not change the group generated by the transpositions ${t'_i}$. 
One has $\Sigma_1=\cup_{i=1}^{q}\Delta_{i}$, so $d=|s|+q$. Adding 
$2|s|$ to both sides of the inequality $n-|s|\geq 2(q-1)$ we obtain 
the equivalent inequality $n+|s|\geq 2d-2 = 2|\Sigma_1|$.
\end{proof}
\begin{dfn}\label{s4.43}
A transitive subgroup $G\subset S_d$ is called \emph{imprimitive} if 
there is a decomposition $\{1,\ldots ,d\}=\bigsqcup_{i=1}^{k}\Sigma_i$ 
such that for $\forall g\in G$ and $\forall i$ one has 
${\Sigma_i}^{g}=\Sigma_j$ for some $j$ and furthermore $\# \, 
\Sigma_i=e$ for $\forall i$ where $1<e<d$. A transitive subgroup 
$G\subseteq S_d$ which is not imprimitive is called \emph{primitive}. 
\end{dfn}
If 
$d$ is a prime number clearly every transitive subgroup of $S_d$ is primitive.
\begin{lem}\label{s4.43a}
Suppose that $G\subseteq S_d$ is a primitive group which contains a 
transposition. Then $G=S_d$. In particular if $d$ is a prime number 
every transitive subgroup of $S_d$ which contains a transposition is 
equal to $S_d$.
\end{lem}
\begin{proof}
Let us prove that if $G$ contains a nontrivial symmetric subgroup 
$S(\Sigma)$ with $\# \, \Sigma <d$ then $G$ contains a $S(\Sigma')$ 
with $\# \, \Sigma'=\# \, \Sigma +1$. Indeed, since $G$ is primitive 
there are $g,h\in G$ such that $0<\# (\Sigma^{g}\cap \Sigma^{h})<\# \, 
\Sigma$. Hence if $f=hg^{-1}$ one has $0<\#(\Sigma \cap \Sigma^{f})<\# 
\, \Sigma$. Let $x^{f}\in \Sigma \cap \Sigma^{f},\; y^{f}\in \Sigma^{f} - 
 \Sigma$. Then $x,y \in \Sigma,\; x^{f}\in \Sigma ,\; 
y^{f}\notin \Sigma$. So $\Sigma'=\Sigma \cup \{x,y\}\supsetneqq 
\Sigma$ and $S(\Sigma')=\langle S(\Sigma),f^{-1}(xy)f\rangle$ is contained in 
$G$. We conclude the proof of the lemma by induction on $\# \, \Sigma$ 
starting from $\Sigma = \{a,b\}$ with $(ab)\in G$.
\end{proof}
To each connected component of $\mathcal{H}_{d,n}(Y,b_0)$ is 
associated a fixed monodromy group $G\subseteq S_d$ and similarly to 
each connected component of $\mathcal{H}_{d,n}(Y)$ is associated a 
conjugacy class of subgroups.
\begin{dfn}\label{s4.44}
Let $d\geq 2,\; n>0$. We denote by $\mathcal{H}^0_{d,n}(Y,b_0)$ and 
$\mathcal{H}^0_{d,n}(Y))$ the union of connected components of
$\mathcal{H}^{(n)}_{d}(Y,b_0)$ and $\mathcal{H}^{(n)}_{d}(Y)$ respectively
 which correspond to coverings with primitive 
monodromy groups (therefore equal to $S_d$ by Lemma~\ref{s4.43a}). 
\end{dfn}
The following theorem is due in the case $n\geq 2d$ to Graber, Harris 
and Starr \cite{GHS}.
\begin{thm}\label{s4.45}
Let $Y$ be a smooth, projective curve. Let $b_0\in Y$. If $n\geq 2d-2$ 
the Hurwitz spaces $\mathcal{H}^0_{d,n}(Y,b_0)$ and 
$\mathcal{H}^0_{d,n}(Y)$ are irreducible.
\end{thm}
\begin{proof}
The case $Y\cong \mathbb{P}^1$ is classical and due to Clebsch and Hurwitz 
\cite{Cl}, \cite{Hu}. Modern proofs may be found in \cite{Mo} p.368 
and \cite{Vo} p.197. Suppose that $g(Y)\geq 1$. Let us first consider  
$\mathcal{H}^0_{d,n}(Y,b_0)$. Since $\mathcal{H}^0_{d,n}(Y,b_0)$ is 
smooth in order to prove its irreducibility it suffices to prove it is
connected. By Theorem~\ref{s2.30c} it suffices to show that every 
Hurwitz system 
$(t_1,\ldots,t_n;\lambda_1,\mu_1,\ldots,\lambda_g,\mu_g)$ with 
monodromy group $S_d$ is braid-equivalent to 
\begin{equation}\label{es4.45a}
((12),(12),(13),(13),\ldots,(1d),(1d);1,1,\ldots,1,1)
\end{equation}
where each $(1i)$ appears twice if $2\leq i\leq d-1$ and $(1d)$ 
appears $n-2d+4$ times.
\begin{stepb}\label{b1}
We claim that every Hurwitz system 
$(t_1,\ldots,t_n;\lambda_1,\mu_1,\ldots,\lambda_g,\mu_g)$ is 
braid-equivalent to a Hurwitz system 
$(\tilde{t}_1,\ldots,\tilde{t}_n;\lambda_1,\mu_1,\ldots,\lambda_g,\mu_
g)$ such that $\langle \tilde{t}_1,\ldots,\tilde{t}_n\rangle$ is a 
transitive group. Indeed, let $s=[\lambda_{1},\mu_{1}]\cdots 
[\lambda_{g},\mu_{g}]$ and let 
$\underline{\Sigma}=\{\Sigma_1,\ldots,\Sigma_m\}$ be the domains of 
transitivity of the group $\langle t_{1},\ldots ,t_{n} \rangle$. One 
has $|\underline{\Sigma}|=d-m$, so if $\langle t_{1},\ldots ,t_{n} 
\rangle$ is not transitive then $|\underline{\Sigma}|<d-1$. In this 
case the inequality $n+|s|>2|\underline{\Sigma}|$ is satisfied, so 
according to Lemma~\ref{s4.42} one may replace $(t_1,\ldots,t_n)$ by a 
braid-equivalent sequence $(t'_1,\ldots,t'_n)$ such that 
$t'_{n-1}=t'_n=(ab)$ and $\langle t'_{1},\ldots ,t'_{n-2} \rangle = 
\langle t_{1},\ldots ,t_{n} \rangle$, so $H=\langle 
t'_1,\ldots,t'_{n-2};\lambda_1,\mu_1,\ldots,\lambda_g,\mu_g\rangle = 
S_d$. Suppose that $a,b\in \Sigma_i$. Let $h\in H$ be such that $a^{h}\in 
\Sigma_1,\; b^{h}\notin \Sigma_1$. Then according to Main 
Lemma~\ref{s3.32} the Hurwitz system 
$(t'_1,\ldots,t'_{n-2},t'_{n-1},t'_n;
\lambda_1,\mu_1,\ldots,\lambda_g,\mu_g)$ is braid-equivalent to 
$(t'_1,\ldots,t'_{n-2},(t'_{n-1})^{h},(t'_n)^{h};
\lambda_1,\mu_1,\ldots,\lambda_g,\mu_g)$. If 
$\underline{\Sigma}'=\{\Sigma'_1,\ldots,\Sigma'_r\}$ are the domains 
of transitivity of $\langle 
t'_1,\ldots,t'_{n-2},(t'_{n-1})^{h},(t'_n)^{h}
\rangle$ then clearly $|\underline{\Sigma}'|>|\underline{\Sigma}|$. 
Repeating this argument after a finite number of steps one obtains a 
braid-equivalent Hurwitz system 
$(\tilde{t}_1,\ldots,\tilde{t}_n;\lambda_1,\mu_1,\ldots,\lambda_g,\mu_
g)$ with transitive $\langle \tilde{t}_1,\ldots,\tilde{t}_n\rangle$.
\end{stepb}
\begin{stepb}\label{b2}
Let $(t_1,\ldots,t_n;\lambda_1,\mu_1,\ldots,\lambda_g,\mu_g)$ be a 
Hurwitz system. Suppose that 
$\sum_{\ell =1}^{g}(|\lambda_{\ell}|+|\mu_{\ell}|)>0$. We claim there 
is a braid-equivalent sequence 
$(t'_1,\ldots,t'_n;\lambda'_1,\mu'_1,\ldots,\lambda'_g,\mu'_g)$ such 
that 
$\sum_{\ell =1}^{g}(|\lambda'_{\ell}|+|\mu'_{\ell}|)<
\sum_{\ell =1}^{g}(|\lambda_{\ell}|+|\mu_{\ell}|)
$. Let us first suppose  that $\lambda_{1}\neq 1$. Decomposing $\lambda_{1}$ 
into a product of nontrivial independent cycles and choosing one of 
them $(ab\ldots c)$ we have $|(ab)\lambda_{1}|<|\lambda_{1}|$. 
According to Step~\ref{b1} we may replace 
$(t_1,\ldots,t_n;\lambda_1,\mu_1,\ldots,\lambda_g,\mu_g)$ by a 
braid-equivalent system 
$(\tilde{t}_1,\ldots,\tilde{t}_n;\lambda_1,\mu_1,\ldots,\lambda_g,\mu_
g)$ such that $\tilde{t}_1,\ldots,\tilde{t}_n$ generate a transitive 
group. Then according to Mochizuki's Lemma~\ref{s4.37a} we may replace 
$(\tilde{t}_1,\ldots,\tilde{t}_n)$ by a braid-equivalent sequence 
$((ab),t'_2,\ldots,t'_n)$. Applying the braid move $\tau_{11}''$ (cf. 
Eq.\eqref{es2.30b}) we transform $\lambda_{1}$ into $(ab)\lambda_{1}$. 
Suppose that $\lambda_{1}=1,\; \mu_{1}\neq 1$. We reason in the same way 
and use the braid move $\rho_{11}'$ which transforms $((ab),\ldots 
;1,\mu_{1},\ldots )$ into $(t'_1,\ldots ;1,(ab)\mu_{1},\ldots )$ 
according to Corollary~\ref{s2.30a}~(i). Similarly if $\lambda_{k}\neq 
1$ and $\lambda_{\ell}=\mu_{\ell}=1$ for $\forall \ell \leq k-1$ we 
have $u_{k-1}=1$ and applying the braid move $\tau_{1k}''$ we transform 
$((ab),\ldots ;\ldots ,\lambda_{k},\mu_{k},\ldots )$ into 
$(t'_1,\ldots;\ldots,(ab)\lambda_{k},\mu_{k},\ldots)$ thus decreasing 
$|\lambda_{k}|$. If $\mu_{k}\neq 1$, 
$\lambda_{\ell}=\mu_{\ell}=1$ for $\forall \ell \leq k-1$ and 
$\lambda_{k}=1$ we apply $\rho_{1k}'$ to transform 
$((ab),\ldots;\ldots,1,\mu_{k},\ldots)$ into 
$(t_1',\ldots;\ldots,1,(ab)\mu_{k},\ldots)$ thus decreasing 
$|\mu_{k}|$.
\end{stepb}
\begin{stepb}
Starting with 
$(t_1,\ldots,t_n;\lambda_1,\mu_1,\ldots,\lambda_g,\mu_g)$ and applying 
several times Step~\ref{b2} one obtains a braid-equivalent Hurwitz 
system $(\tilde{t}_1,\ldots,\tilde{t}_n;1,1,\ldots,1,1)$. Here 
$\tilde{t}_{1}\cdots \tilde{t}_{n}=1$, so applying the argument in the 
case $g(Y)=0$ (see e.g. \cite{Mo} p.368, or \cite{Vo} p.197) one 
obtains the initial Hurwitz system is braid-equivalent to 
Eq.\eqref{es4.45a}. 
\end{stepb}
The connectedness of both $\mathcal{H}^0_{d,n}(Y,b_0)$ and 
$\mathcal{H}^0_{d,n}(Y)$ follows now from Theorem~\ref{s2.30c}.
\end{proof}
Theorem~\ref{s4.45} may be generalized in a straightforward manner to 
coverings $\pi :X\to Y$ which have simple ramifications in all except 
possibly one discriminant point. The conjugacy classes of $S_d$ 
correspond bijectively to partitions of $d$, namely 
$\underline{e}=\{e_1,e_2,\ldots,e_q\}$ where $e_1\geq e_2\geq \ldots 
\geq e_q\geq 1$ and $e_1+\cdots +e_q=d$. To each partition one 
associates the orbit of the permutation 
\begin{equation}\label{es4.49}
\epsilon = (1\, 2\ldots e_1)(e_1+1\ldots e_1+e_2)\ldots 
(\sum_{i=1}^{q-1}e_i\:+1\ldots d)
\end{equation}
Given a partition $\underline{e}$ let us denote by 
$\mathcal{H}_{d,n,\underline{e}}(Y,b_0)$ and 
$\mathcal{H}_{d,n,\underline{e}}(Y)$) the Hurwitz spaces which 
parametrize equivalence classes $[X\to Y,\: \phi]$ and $[X\to 
Y]$ respectively, where $X$ is connected and $\pi :X\to Y$ is a covering of degree 
$d$ with $n$ discriminant points of simple ramification and one 
additional discriminant point whose local monodromy belongs to the 
conjugacy class of Eq.\eqref{es4.49}. Let 
$\mathcal{H}^0_{d,n,\underline{e}}(Y,b_0)$
and 
$\mathcal{H}^0_{d,n,\underline{e}}(Y)$ 
parametrize such 
coverings with primitive monodromy group.
\begin{thm}\label{s4.49}
Let $Y$ be a smooth, projective curve. Let $b_0\in Y$.
Suppose that $n\geq 2d-2$. Let $\underline{e}$ be an 
arbitrary partition of $d$. Then the Hurwitz spaces 
$\mathcal{H}^0_{d,n,\underline{e}}(Y,b_0)$ and 
$\mathcal{H}^0_{d,n,\underline{e}}(Y)$  are irreducible.
\end{thm}
\begin{proof}
If $\underline{e}$ is the trivial partition (that is $e_i=1$ for $\forall i$) 
this is the 
content of Theorem~\ref{s4.45}. 
Suppose that $e_1\geq 2$. Let us denote 
the permutation Eq.\eqref{es4.49} by 
$\epsilon = (1_12_1\ldots (e_1)_1)\linebreak
(1_2\ldots (e_2)_2)\ldots$. The 
theorem will be proved if we can show that every Hurwitz system 
$(t_1,\ldots,t_n,t_{n+1};\lambda_1,\mu_1,\ldots,\lambda_g,\mu_g)$ of 
the considered type may be reduced after a finite number of braid 
moves of the types 
$\sigma_{j}',\sigma_{j}'',\rho_{ik}',\rho_{ik}'',\tau_{ik}',
\tau_{ik}''$ to the normal form
\[
(\mathcal{Z},t_{N+1}',\ldots ,t'_n,\epsilon^{-1};1,1,\ldots ,1,1)
\]
where the $n$-tuple of transpositions $(\mathcal{Z},t_{N+1}',\ldots 
,t'_n)$ is the one defined in \ref{s4.40a} and 
Theorem~\ref{s4.40}. First using braid moves of type $\sigma_{j}'$ 
we may replace the original Hurwitz system by one for which $t_{n+1}$ 
belongs to the orbit of $\epsilon$. We then let 
$s=[\lambda_{1},\mu_{1}]\cdots [\lambda_{g},\mu_{g}]t^{-1}_{n+1}$ and 
repeating the arguments of the proof of Theorem~\ref{s4.45} we obtain 
a braid-equivalent Hurwitz system 
$(\tilde{t}_1,\ldots,\tilde{t}_n,\tilde{t}_{n+1};1,1,\ldots,1,1)$ with 
$\tilde{t}_{n+1}=t_{n+1}$. We are further allowed to apply only braid 
moves of types $\sigma_{j}'$ and $\sigma_{j}''$. The next step is to 
replace the obtained Hurwitz system by one in which at the $(n+1)$-th 
place stays $\epsilon^{-1}$. In fact the permutation $\epsilon^{-1}$ 
has the same cyclic type as $\tilde{t}_{n+1}$, so $\epsilon^{-1}=
a^{-1}\tilde{t}_{n+1}a$. Let $a=\tau_{1}\cdots \tau_{r}$ where 
$\tau_{i}$ are transpositions. By the hypothesis and Lemma~\ref{s4.43a}
we have $\langle \tilde{t}_1,\ldots,\tilde{t}_n\rangle=S_d$. So using 
Mochizuki's Lemma~\ref{s4.37a} we may replace 
$(\tilde{t}_1,\ldots,\tilde{t}_n)$ by a braid-equivalent $n$-tuple 
$(\ldots,\tau_{1})$. We then have $(\tau_{1},\tilde{t}_{n+1})\sim 
(\tau_{1}\tilde{t}_{n+1}\tau_{1},\tau_{1})\sim 
(\tau_{1}',\tau_{1}\tilde{t}_{n+1}\tau_{1})$, so 
$(\tilde{t}_1,\ldots,\tilde{t}_n,\tilde{t}_{n+1})\sim 
(\ldots,\tau_{1}\tilde{t}_{n+1}\tau_{1})$. Repeating this argument 
with $\tau_{2},\tau_{3},\ldots,\tau_{r}$ we obtain a braid-equivalent 
$(n+1)$-tuple whose $(n+1)$-th element is $\epsilon^{-1}$ and the 
product of the first $n$ equals $\epsilon$. Using a finite number of 
braid moves among the first $n$ transpositions we obtain the required 
normal form according to Theorem~\ref{s4.40}. 
\end{proof}
So far we worked with coverings with primitive monodromy groups. We now 
wish to consider the imprimitive case.
\begin{lem}\label{s4.52}
Let $G\subset S_d$ be a transitive imprimitive  subgroup which 
contains a transposition. Then there is a unique decomposition 
$\{1,\ldots ,d\}=\Sigma_1\sqcup \ldots \sqcup \Sigma_k$ as in 
Definition~\ref{s4.43} so that $G_i=G\cap S(\Sigma_i)$ is a primitive 
subgroup of $S(\Sigma_i)$ for $\forall i$. Furthermore the following 
properties hold.
\renewcommand{\theenumi}{\roman{enumi}}
\begin{enumerate}
\item All transpositions contained in $G$ form one conjugacy $G$-orbit 
$T$.
\item $G_1\cdot G_2\cdots G_k\ =\ G_1\times G_2\times \cdots \times 
G_k = \langle T\rangle$.
\item $G_i = S(\Sigma_i)$ for $\forall i$.
\end{enumerate}
\end{lem}
\begin{proof}
Let $(ab)\in G$ and let $T=(ab)^{G}$. Let $H=\langle T\rangle$ and let 
$\Sigma_1,\ldots,\Sigma_k$ be the domains of transitivity of $H$. 
Since $H$ is a normal subgroup every $g\in G$ permutes the orbits 
$\{\Sigma_1,\ldots,\Sigma_k\}$. Furthermore $H=S(\Sigma_1)\times 
\cdots \times S(\Sigma_k)$, so $G_i=G\cap S(\Sigma_i)=S(\Sigma_i)$. 
All transpositions which belong to $H$ form one $G$-orbit since this 
property holds for $S(\Sigma_i)$. Let $(\alpha \, \beta)$ be an 
arbitrary transposition of $G$. We claim it is impossible that 
$\alpha\in \Sigma_i,\; \beta \in \Sigma_j$ for $i\neq j$. Indeed, if 
this were the case letting $g=(\alpha \, \beta)$ we would have 
${\Sigma_i}^{g}\cap \Sigma_j\neq \emptyset$, so 
${\Sigma_i}^{g}=\Sigma_j$. If $\alpha'\in \Sigma_i, \alpha'\neq 
\alpha$, then ${\alpha'}^{g}=\alpha'\in \Sigma_i$ and on the other hand 
${\alpha'}^{g}\in \Sigma_j$, so $\Sigma_i\cap \Sigma_j\neq \emptyset$ 
which is absurd. We obtain $\{\alpha\, \beta\}\subseteq \Sigma_i$ for 
some $i$, so $(\alpha\, \beta)\in T=(ab)^{G}$. This proves (i). Suppose now that $\{1,\ldots 
,d\}=\Sigma'_1\sqcup \ldots \sqcup \Sigma'_{\ell}$ is an arbitrary 
decomposition as required in Definition~\ref{s4.43}. Let $(\alpha\, 
\beta)$ be a transposition in $G$. Then the argument above shows 
that $\{\alpha,\beta\}\subseteq \Sigma'_i$ for some $i$. This shows 
that
$T=(\alpha\, \beta)^{G}\subseteq G' = G'_1\times \cdots \times 
G'_{\ell}$ where $G'_j=G\cap S(\Sigma'_j)$. Thus each orbit $\Sigma_i$ 
of $H=\langle T\rangle$ is contained in some orbit $\Sigma'_j$. 
Assuming $G'_j$ is primitive subgroup of $S(\Sigma'_j)$ we conclude that
$\Sigma_i=\Sigma'_j$. This proves the uniqueness statement. The other 
properties were already proved.
\end{proof}
\begin{cor}\label{s4.54}
Let $\pi :X\to Y$ be a covering of smooth, irreducible, projective 
curves of degree $d$. Suppose that at least one of the discriminant points 
is simple and suppose that the monodromy group of the covering is $\neq 
S_d$. Then there exists a unique decomposition 
$X \overset{\pi_1}{\lto}\tilde{Y}\overset{\pi_2}{\lto}Y,\; \pi = 
\pi_2\circ \pi_1$, where $d_i=\deg \pi_i > 1$ for $i=1,2$ and 
$\pi_1:X\to \tilde{Y}$ has primitive monodromy group. If 
$\pi : X\to Y$ 
has $n$ simple discriminant points so does $\pi_1:X\to \tilde{Y}$ and 
$\pi_2$ induces a bijective correspondence between these two sets. 
\end{cor}
\begin{proof}
One chooses an unramified point $b_0\in Y$, and a bijection $\phi: 
\pi^{-1}(b_0)\overset{~}{\lto}\{1,\ldots ,d\}$. One applies 
Lemma~\ref{s4.52} to the monodromy group of $X\to Y$ and obtains a decomposition $X\to \tilde{Y}\to Y$. 
Replacing $\phi$ by another bijection results in replacing the 
monodromy group $G$ by a conjugate one $G'=s^{-1}Gs$. The conjugation 
by $s\in S_d$ transforms the $G$-orbit of transpositions contained in 
$G$ to the $G'$-orbit of transpositions contained in $G'$. Therefore 
$s$ transforms the corresponding decompositions of $\{1,\ldots ,d\}$:\quad 
$\Sigma'_i = {\Sigma_i}^{s}$. This proves the uniqueness 
of the decomposition $X\to \tilde{Y}\to Y$ with the required 
properties. The last statement of the corollary is obvious.
\end{proof}
\begin{thm}\label{s4.55}
Let $Y$ be a smooth, projective curve of genus $\geq 1$. 
Let $d\geq 2,\; n>0,\; n\equiv 0(mod\; 2)$. If $d$ is prime then 
$\mathcal{H}_{d,n}(Y) = \mathcal{H}^0_{d,n}(Y)$. If $d$ is not prime
let $d'$ be its maximal divisor  $\neq d$. 
 Suppose that $n\geq 2d'-2$. Then the connected 
components of $\mathcal{H}_{d,n}(Y)$ which correspond to simple coverings 
with monodromy groups $\neq S_d$ are in one-to-one correspondence with 
the equivalence classes of unramified coverings $[\tilde{Y}\to Y]$ 
of degrees $d_2|d$ where $d_2\neq 1, d$.
\end{thm}
\begin{proof}
If $d$ is prime every transitive subgroup 
of $S_d$ is primitive, so 
$\mathcal{H}_{d,n}(Y) = \mathcal{H}^0_{d,n}(Y)$. Suppose that $d$ is not 
prime.
Let $\pi :X\to Y$ be a simple covering of degree $d$, simply ramified 
in $n$ points, with monodromy group $\neq S_d$.
 Then the decomposition 
$X \overset{\pi_1}{\lto}\tilde{Y}\overset{\pi_2}{\lto}Y$ of 
Corollary~\ref{s4.54} is with \'{e}tale $\pi_2$. Furthermore by the 
uniqueness of this decomposition if $f:X\to X'$ defines an equivalence 
of coverings of $Y$ then $f$ is induced by an equivalence 
$g:\tilde{Y}\to \tilde{Y}'$. Given an \'{e}tale covering 
$p_2:\tilde{Y}\to Y$ of degree $d_2\neq 1,d$ we consider the Hurwitz 
space $\mathcal{H}^0_{d_1,n}(\tilde{Y})$ where $d_1=
d/d_2$. It is 
connected by Theorem~\ref{s4.45}. We then consider the Zariski open 
subset consisting of $[X\to \tilde{Y}]$ with discriminant points 
belonging to different fibers of $\pi_2 :\tilde{Y}\to Y$. Every such 
$X\to \tilde{Y}$ yields a simple covering of $Y$. In this way one 
obtains a connected component of $\mathcal{H}_{d,n}(Y)$ whose 
conjugacy class of monodromy groups is imprimitive. Vice versa by 
Corollary~\ref{s4.54} every connected component of 
$\mathcal{H}_{d,n}(Y)$ whose associated conjugacy class of monodromy 
groups is imprimitive is obtained in this way.
\end{proof}
\begin{exa}\label{s4.57}
Let $d=4$ and let $Y$ be a curve of genus $g\geq 1$. The connected 
unramified coverings of $Y$ of degree 2 are classified up to 
equivalence by the points of order 2 in the Jacobian $J(Y)$. Applying 
Theorem~\ref{s4.55} we obtain if $n\geq 2$ there are 
$2^{2g}-1$ different connected components of 
$\mathcal{H}_{4,n}(Y)$ which parametrize coverings with imprimitive 
monodromy group. The monodromy group in this case is isomorphic to the 
dihedral group $D_4$
\end{exa}

\section{The case $n=2d-4,\quad g\geq 1$}\label{s5}
Unless otherwise specified in this section we shall work with Hurwitz 
systems $(t_1,\ldots,t_n;\linebreak
\lambda_1,\mu_1,\ldots,\lambda_g,\mu_g)$ 
where $n\geq 2,\; g\geq 1$ and $t_1,\ldots,t_n$ are transpositions in 
$S_d$. 
Let $T=(t_1,\ldots,t_n)$ be a sequence of transpositions in $S_d$ and let 
$s=t_{1}\cdots t_{n}$. Suppose that the subgroup 
$\langle T\rangle=\langle t_1,\ldots,t_n\rangle$ 
has only one 
nontrivial domain of transitivity 
$\Sigma$ and let $e=\#\, \Sigma$. We have the inequality $n\geq e-1,\; 
n+|s|\geq 2(e-1)$ (cf. Lemma~\ref{s4.42}) so if $e=\#\, \Sigma$ is 
fixed the minimum for $n$ is reached for $n=e-1$, and then $|s|=e-1$. This 
happens if and only if $s=t_{1}\cdots t_{n}$ is a cycle of order 
$n+1$. 
\begin{dfn}\label{s5.58bis}
We call a sequence of transpositions $T=(t_1,\ldots,t_n)$ 
\emph{minimal} if 
$s=t_{1}\cdots t_{n}$ is a cycle of order $n+1$. 
\end{dfn}
\begin{lem}\label{s5.59bis}
Let $T=(t_1,\ldots,t_n)$ be a sequence of transpositions of $S_d$ such 
that $\langle T\rangle=\langle t_1,\ldots,t_n\rangle$ has a unique 
nontrivial domain of transitivity $\Sigma$. Suppose that $T$ is 
not minimal. Then for any $a\in \Sigma$ there exists a 
braid-equivalent sequence $T'=(t'_1,\ldots,t'_n)$ such that 
$t'_{n-1}=t'_n=(ab)$ for some $b\in \Sigma,\; b\neq a$. 
\end{lem}
\begin{proof}
Without loss of generality we may assume that $\Sigma=\{1,\ldots ,d\}$. Let 
$s=t_{1}\cdots t_{n}$ and let $\Delta_{1},\ldots,\Delta_{q}$ be the 
domains of transitivity of $s$. Replacing $(t_1,\ldots,t_n)$ by a 
braid-equivalent sequence we may assume it has the normal form of 
Theorem~\ref{s4.40}. The non-minimality hypothesis $n>|s|$ implies 
that $T$ contains $(n-|s|)/2$ pairs of the form $(1_12_1)$ if 
$q=1$ or $(1_11_i)$ with $2\leq i\leq q$ if $q\geq 2$. Consider the 
cyclic group $\langle s\rangle$. If $u\in \langle s\rangle$ then 
${t_{1}}^{u}\cdots {t_{n}}^{u}=s$, so the sequence 
$({t_1}^{u},\ldots,{t_n}^{u})$ is braid-equivalent to 
$(t_1,\ldots,t_n)$ (cf. \cite{Kl}, \cite{Mo}). If $q=1$ we may find 
$u\in \langle s\rangle$ such that $(1_1)^{u}=a$. Here we let 
$b=(2_1)^{u}$. If $q\geq 2$ and $a\in \Delta_{1}$ we find $u$ such 
that $(1_1)^{u}=a$ and we let $b=(1_2)^{u}$. If $a\in \Delta_{i}$ with 
$i\geq 2$ we find $u\in \langle s\rangle$ such that $(1_i)^{u}=a$ and we
let $b=(1_1)^{u}$. In each case the sequence 
$({t_1}^{u},\ldots,{t_n}^{u})$ contains the adjacent pair $(ab),(ab)$. 
Moving it to the end we obtain the sequence required in the lemma. 
\end{proof}
\begin{block}\label{s5.60bis}
Given a Hurwitz system 
$(t_1,\ldots,t_n;\lambda_1,\mu_1,\ldots,\lambda_g,\mu_g)$ where 
$t_i,\lambda_{k},\mu_{k}\in S_d$ let 
$\underline{\Sigma}=\{\Sigma_1,\Sigma_2,\ldots,\Sigma_m\}$ be the 
domains of transitivity of $G_1=\langle t_{1},\ldots ,t_{n} \rangle$ 
ordered in such a way that 
$\# \, \Sigma_1\geq \# \, \Sigma_2\geq 
\cdots \geq \# \, \Sigma_m$. Let 
\begin{equation}\label{es5.60bis}
\ell(\underline{\Sigma})\ =\ (\# \, \Sigma_1,\# \, \Sigma_2,\ldots,\# 
\, \Sigma_m,0,0,\ldots)
\end{equation}
be the associated partition of $d$. Given $\Sigma_i$ let $T_i$ be the 
subsequence of $T=(t_1,\ldots,t_n)$ composed of those transpositions 
which move points of $\Sigma_i$. Let $n_i= \# \, T_i$. We call 
$\Sigma_i$ minimal or non-minimal if the sequence $T_i$ is 
minimal or non-minimal respectively. The trivial case $\# \, 
\Sigma_i=1,\; T_i=\emptyset$ is assumed minimal.
\end{block}
\begin{lem}\label{s5.61bis}
Let $(t_1,\ldots,t_n;\lambda_1,\mu_1,\ldots,\lambda_g,\mu_g)$ be a 
Hurwitz system with monodromy group $S_d$. Let 
$s=[\lambda_{1},\mu_{1}]\cdots [\lambda_{g},\mu_{g}]$ and let 
$G_1=\langle t_{1},\ldots ,t_{n} \rangle \subseteq S_d$. Let 
$\underline{\Sigma}=\{\Sigma_1,\ldots,\Sigma_m\}$ be the domains of 
transitivity of $G_1$ ordered as in \ref{s5.60bis}. Suppose that
$\ell(\underline{\Sigma})$ is maximal in lexicographic order among the 
partitions associated with all Hurwitz systems of the type 
$(t'_1,\ldots,t'_n;\lambda_1,\mu_1,\ldots,\lambda_g,\mu_g)$ 
braid-equivalent to the given one. Then 
$\underline{\Sigma}=\{\Sigma_1\}$ (equivalently $\langle t_{1},\ldots ,t_{n} 
\rangle$ is transitive) if and only if $n+|s|\geq 2d-2$. If 
$n+|s|<2d-2$ then $n+|s|=2|\underline{\Sigma}|=2(d-m)$. Furthermore 
$n_1+|s|_{\Sigma_1}|=2(\# \, \Sigma_1 -1)$ and for each $i\geq 2$ the 
$G_1$-orbit
$\Sigma_i$ is minimal, equivalently $n_i=|s|_{\Sigma_i}|=\# \, \Sigma_i -1$. 
  In particular if $s=1$ then $\langle t_{1},\ldots ,t_{n} \rangle$ has a 
single nontrivial domain of transitivity.
\end{lem}
\begin{proof}
We proved in Lemma~\ref{s4.42} that if $\langle t_{1},\ldots ,t_{n} 
\rangle$ is transitive then $n+|s|\geq 2d-2$.  That the maximality 
of $\ell(\underline{\Sigma})$ and $n+|s|\geq 2d-2$ imply that
$\underline{\Sigma}=\{\Sigma_1\}$ is a fact evident from the proof of 
Theorem~\ref{s4.45} (cf. Step~\ref{b1}). We shall further assume that 
$n+|s|<2d-2$, so $G_1$ is not transitive. Suppose, by way of 
contradiction, that there is a non-minimal $T_k$ with $k\geq 2$. Let 
$\Sigma'=\Sigma_{k+1}\cup\cdots \cup \Sigma_{m}$, so $\{1,\ldots 
,d\}=\Sigma_1\cup \cdots \cup \Sigma_k \cup \Sigma'$. Let $H_1=\langle 
T_k\rangle,\: H_2=\langle 
T_1,\ldots,\widehat{T_k},\ldots,T_{m},
\lambda_1,\mu_1,\ldots,\lambda_g,\mu_g \rangle$, so $\langle 
H_1,H_2\rangle=G=S_d$. For every $b\in \Sigma_k$ one has 
$b^{G}=\{1,\ldots ,d\}$. Let $a\in \Sigma_k$ and $h\in G=\langle 
H_1,H_2\rangle$ be chosen in such a way
 that $a^{h}\in 
\Sigma_1\cup \cdots \cup \Sigma_{k-1}$ and furthermore  the 
length of the factorization 
$h=h_1h_2\cdots h_r$ with $h_i\in H_1$ or $H_2$ is minimal possible 
among all $b\in \Sigma_k$ and all $g\in G$ satisfying $b^{g}\in 
\Sigma_1\cup \cdots \cup \Sigma_{k-1}$. We claim $r=1$ and $h_1\in 
H_2$. In fact by minimality we have $h_1\in H_2,\: h_r\in H_2$ and for 
each $i=1,\ldots,r-1$ the adjacent $h_i,h_{i+1}$ do not belong to the 
same subgroup $H_1$ or $H_2$. Suppose that $r\geq 2$. Then by minimality 
$a^{h_1\cdots h_{r-1}}\in \Sigma_{k}\cup \Sigma'$. If 
$a^{h_1\cdots h_{r-1}}\in \Sigma_k$ we might replace $a$ by 
$a'=a^{h_1\cdots h_{r-1}}$, so $a^{h_1\cdots h_{r-1}}\in \Sigma'$. The 
group $H_1$ acts trivially on $\Sigma'$, therefore since $h_{r-1}\in 
H_1$ we have $a^{h_1\cdots h_{r-1}}=a^{h_1\cdots h_{r-2}}$. Hence 
$a^h=a^{h_1\cdots (h_{r-2}h_{r})}=a^{h'}$ where $h'$ has a shorter 
factorization then $h$. This is an absurd, so $r=1$. We have thus 
proved there exists an $a\in \Sigma_k$ and $h\in H_2$ such that 
$a^{h}\in \Sigma_i$ for some $i<k$. According to Lemma~\ref{s5.59bis} 
we may replace $T_k$ by a braid-equivalent sequence which contains the 
adjacent pair $(ab),(ab)$ for some $b\in \Sigma_k$. By 
Main Lemma~\ref{s3.32} replacing $(ab),(ab)$ by $(ab)^{h},(ab)^{h}$ one 
obtains a braid-equivalent Hurwitz system. The following cases may 
occur.
\begin{casea}
\emph{$b^{h}\in \Sigma_j$ for some $j\neq i,\; 1\leq j\leq m$}.\quad If 
$i<j$ then the new Hurwitz system has $\Sigma_i\cup \Sigma_j$ as a 
domain of transitivity, so the corresponding partition 
Eq.\eqref{es5.60bis} is greater than $\ell(\underline{\Sigma})$ in 
lexicographical order. If $j<i$ we interchange $a$ and $b$ and make 
the same conclusion. In both cases we obtain a contradiction with the 
maximality of $\ell(\underline{\Sigma})$. 
\end{casea}
\begin{casea}
$b^{h}\in \Sigma_i$.\quad Let $(ab)^{h}=(xy)$. Moving the pair 
$(xy),(xy)$ to the right of $T_i$ we obtain a braid-equivalent sequence 
$\tilde{T}_i=T_i,(xy),(xy)$ with the same domain of transitivity as 
$T_i$. Canceling the pair from the obtained Hurwitz system does not 
change the monodromy group $G$. Let us choose $g\in G=S_d$ such that 
$x^{g}\in \Sigma_i,\; y^{g}\notin \Sigma_i$. Applying Main
Lemma~\ref{s3.32} we obtain a new Hurwitz system 
$(\tilde{t}_1,\ldots,\tilde{t}_n;\lambda_1,\mu_1,\ldots,\lambda_g,\mu_
g)$ whose associated partition $\ell(\underline{\Sigma})$ has the 
property that its first $i$ terms are greater in lexicographical order 
then the first $i$ terms of $\ell(\underline{\Sigma})$. This 
contradicts the maximality of the partition 
$\ell(\underline{\Sigma})$. 
\end{casea}
The claim that $T_i$ is minimal for every $i\geq 2$ is proved. According to 
Lemma~\ref{s4.42} we have 
$n_1+|s|_{\Sigma_1}|\geq 2(\# \, \Sigma_1 -1)$ and if this inequality 
is strict we can apply the argument of 
Step~\ref{b1} of Theorem~\ref{s4.45} in order to increase 
$\ell(\underline{\Sigma})$ which is impossible . Thus 
$n_1+|s|_{\Sigma_1}|= 2(\# \, \Sigma_1 -1)$. Finally if $s=1$ then 
$s|_{\Sigma_i}=1$, so minimality of $\Sigma_i$ means $\# \, 
\Sigma_i=1$. Therefore the maximality of $\ell(\underline{\Sigma})$ 
implies that only $\Sigma_1$ may be nontrivial domain of transitivity.
\end{proof}
\begin{block}\label{s5.63bbis}
Given a Hurwitz system 
$(t_1,\ldots,t_n;\lambda_1,\mu_1,\ldots,\lambda_g,\mu_g)$ we let 
$\underline{\lambda}=(\lambda_1,\ldots,\lambda_g)$, 
$\underline{\mu}=(\mu_1,\ldots,\mu_g)$, 
$|\underline{\lambda}|=(|\lambda_1|,\ldots,|\lambda_g|)$, 
$|\underline{\mu}|=(|\mu_1|,\ldots,|\mu_g|)$ (cf. \ref{s4.41a}). We 
consider the class $\mathcal{C}$ of Hurwitz systems braid-equivalent to the 
given one. When dealing with the problem of finding a system in $\mathcal{C}$ of 
a simplest form it makes sense to assume (without loss of generality) that
$(t_1,\ldots,t_n;\lambda_1,\mu_1,\ldots,\lambda_g,\mu_g)$ satisfies 
the following property:

\smallskip
\noindent
(*) \emph{
$(|\underline{\lambda}|,|\underline{\mu}|)$ is minimal in 
lexicographic order among the Hurwitz systems in $\mathcal{C}$ and furthermore 
$\ell(\underline{\Sigma})$ (cf. \ref{s5.60bis}) is maximal in 
lexicographic order among the Hurwitz systems in $\mathcal{C}$ of the type 
$(t'_1,\ldots,t'_n;\lambda_1,\mu_1,\ldots,\lambda_g,\mu_g)$
}
\end{block}
\begin{thm}\label{s5.64bis}
Let $d\geq 3,\; n=2d-4$. Let $Y$ be a smooth, projective curve of 
genus $g\geq 1$. Then the Hurwitz space $\mathcal{H}^0_{d,n}(Y)$ 
parametrizing simple coverings with monodromy groups equal to  
$S_d$ is irreducible. 
\end{thm}
\begin{proof}
By Theorem~\ref{s2.30c} it suffices to show that the equivalence
 class 
(modulo inner automorphisms) of every Hurwitz system 
$(t_1,\ldots,t_n;\lambda_1,\mu_1,\ldots,\lambda_g,\mu_g)$ with 
monodromy group $S_d$ may be reduced by a finite number of braid moves 
to the normal form
\begin{equation}\label{es5.65bis}
[(12),(12),(13),(13),\ldots,(1\, d-1)(1\, d-1);1,1,\ldots,1,(1d)]
\end{equation}
where each $(1i)$ appears twice for $1\leq i\leq d-1$. Without loss of 
generality we may assume that 
$(t_1,\ldots,t_n;\lambda_1,\mu_1,\ldots,\lambda_g,\mu_g)$ has minimal 
$(|\underline{\lambda}|,|\underline{\mu}|)$ and maximal 
$\ell(\underline{\Sigma})$ as in \ref{s5.63bbis}~(*)
\begin{stepc}\label{c1}
\emph{We claim $\lambda_{1}=\cdots =\lambda_{g}=1.$}\quad Let us first 
suppose that $\lambda_{1}\neq 1$. Let $(ab)$ be a transposition such that 
$|(ab)\lambda_{1}|<|\lambda_{1}|$. Let $s=[\lambda_{1},\mu_{1}]\cdots 
[\lambda_{g},\mu_{g}]$. First suppose  that $s\neq 1$. Since $s\in A_d$ we 
have $|s|\geq 2$, thus $n+|s|\geq 2d-2$. By Lemma~\ref{s5.61bis} 
$\langle t_{1},\ldots ,t_{n} \rangle$ is a transitive group, so by 
Mochizuki's Lemma~\ref{s4.37a} we can replace $(t_1,\ldots,t_n)$ by a 
braid-equivalent sequence $(t'_1,\ldots,t'_n)$ such that $t'_1=(ab)$. 
Applying $\tau''_{11}$ (cf. Eq.\eqref{es2.30b}) one transforms 
$\lambda_{1}$ into $(ab)\lambda_{1}$. This contradicts the minimality 
of $(|\underline{\lambda}|,|\underline{\mu}|)$. Suppose that $s=1$. Then 
according to Lemma~\ref{s5.61bis} the subgroup $\langle t_{1},\ldots 
,t_{n} \rangle$ has two orbits:\quad 
$\underline{\Sigma}=\{\Sigma_1,\Sigma_2\}$ where $\# \, 
\Sigma_1=d-1,\; \# \, \Sigma_2=1$. Moreover replacing if necessary 
$(t_1,\ldots,t_n)$ by a braid-equivalent sequence we may assume that 
$t_i=t_{i+1}$ for every $i\equiv 1(mod\, 2)$. Varying over all Hurwitz 
systems with $\lambda_1,\mu_1,\ldots,\lambda_g,\mu_g$ fixed and 
braid-equivalent to the given one we want to figure out which transpositions 
may appear at the first place. Let $\Sigma_1=\{a_1,\ldots,a_{d-1}\},\; 
\Sigma_d=\{a_d\}$. If $\{\alpha,\beta\}\subseteq \Sigma_1$ then by 
Mochizuki's Lemma~\ref{s4.37a} one may replace the sequence 
$(t_1,\ldots,t_n)$ by a braid-equivalent one 
$((a_{\alpha},a_{\beta}),\ldots )$. Let $G_1=\langle t_{1},\ldots 
,t_{n} \rangle$, $G_2=\langle 
\lambda_1,\mu_1,\ldots,\lambda_g,\mu_g\rangle$. By hypothesis 
$G=\langle G_1,G_2\rangle=S_d$ so $\exists h\in G_2$ such that 
${\Sigma_1}^{h}\not \subseteq \Sigma_1$. Let ${\Sigma_1}^{h}=\{1,\ldots 
,d\}-\{a_c\}$, so ${a_d}^{h}=a_c$ with $c<d$. By Main 
Lemma~\ref{s3.32} the Hurwitz system 
$({t_1}^{h},\ldots,{t_n}^{h};\lambda_1,\mu_1,\ldots,\lambda_g,\mu_g)$ 
is braid-equivalent to the given one. Applying again Mochizuki's 
Lemma~\ref{s4.37a} we see that every $(a_{\alpha}b_{\beta})$ with 
$(\alpha,\beta)\neq (c,d)$ may be placed first in some 
braid-equivalent Hurwitz system. We may vary $h\in G_2$ with the property 
${\Sigma_1}^{h}\not \subseteq \Sigma_1$, or equivalently with the property 
${a_d}^{h}\in \Sigma_1$. If the orbit ${a_d}^{G_2}$ has $\geq 3$ 
elements, then an arbitrary transposition $\tau$ may be placed first 
in some braid-equivalent Hurwitz system 
$(t'_1,\ldots,t'_n;\lambda_1,\mu_1,\ldots,\lambda_g,\mu_g)$. Doing 
this for $\tau = (ab)$ and applying the braid move $\tau''_{11}$ we 
transform 
$\lambda_{1}$ into $(ab)\lambda_{1}$ thus obtaining a contradiction 
with the minimality of $(|\underline{\lambda}|,|\underline{\mu}|)$. It 
remains to consider the case $\# \, {a_d}^{G_2}=2$. In this case 
$G_2\subseteq S(\{a_c,a_d\})\times S(\Sigma_1-\{a_c\})$. If 
$\lambda_{1}\in S(\Sigma_1-\{a_c\})$ or 
$\lambda_{1}=(a_ca_d)\lambda_{1}'$ with $\lambda_{1}'\in 
S(\Sigma_1-\{a_c\})$, $\lambda_{1}'\neq 1$ we might decrease 
$|\lambda_{1}|$ transforming $\lambda_{1}$ into $(ab)\lambda_{1}$ with 
an appropriate $(ab)$ as we saw above and this contradicts the 
minimality of $(|\underline{\lambda}|,|\underline{\mu}|)$. It remains 
to consider the case $\lambda_{1}=(a_ca_d)$. Let $e\in 
\{1,\ldots,d-1\},\: e\neq c$. Replacing the $n$-tuple 
$(t_1,\ldots,t_n)$ by a braid-equivalent one we may assume that it equals 
$((a_ca_e),(a_ca_e),\ldots )$. Since $G_2\subseteq S(\{a_c,a_d\})\times 
S(\Sigma_1-\{a_c\})$ we have $[\lambda_{1},\mu_{1}]=1$. So 
$\rho_{11}'': t_1\mapsto \lambda_{1}t_1\lambda_{1}^{-1}$. Applying 
this braid move we obtain 
$((a_ea_d),(a_ca_e),\ldots;\lambda_{1},\mu_{1}',\ldots)$. We may 
replace $(a_ea_d),(a_ca_e)$ by the braid-equivalent pair 
$(a_ca_d),(a_ea_d)$ and then apply $\tau_{11}''$ transforming 
$\lambda_{1}$ into $(a_ca_d)\cdot (a_ca_d) = 1$. This contradicts 
again the minimality of $(|\underline{\lambda}|,|\underline{\mu}|)$. 
This proves $\lambda_{1}=1$. If $\lambda_{1}=\cdots 
=\lambda_{k-1}=1$, but $\lambda_{k}\neq 1$ we have 
$u_{k-1}=[\lambda_{1},\mu_{1}]\cdots [\lambda_{k-1},\mu_{k-1}]=1$. The 
braid move $\tau_{1k}''$ transforms $\lambda_{k}$ into 
$t_1^{-1}\lambda_{k}$ according to Theorem~\ref{s2.27} so the same 
argument as above may be applied proving that $\lambda_{k}\neq 1$ 
contradicts the minimality of 
$(|\underline{\lambda}|,|\underline{\mu}|)$. Therefore 
$\lambda_{1}=\cdots =\lambda_{g}=1$.
\end{stepc}
\begin{stepc}
\emph{We claim the group $G_1=\langle t_{1},\ldots ,t_{n} \rangle$ has 
orbits $\Sigma_1=\{a_1,\ldots,a_{d-1}\}$, $\Sigma_2=\{a_d\}$ and there 
exists $c\in [1,d-1]$ such that every $\mu_{\ell}$ equals either 1 or 
$(a_ca_d)$}. For every Hurwitz system with $\lambda_{1}=\cdots 
=\lambda_{g}=1$ we have by Theorem~\ref{s2.27}~(b) that $\rho_{1k}'$ 
transforms $\mu_{k}$ into $t_1^{-1}\mu_{k}$. Assume that $\mu_{1}=\cdots 
=\mu_{k-1}=1$ and $\mu_{k}\neq 1$. 
Upon a substitution $\tau \leftrightarrow \rho$
we may use arguments similar to 
those of Step~\ref{c1} trying to decrease $|\mu_{k}|$ and thus 
obtaining a contradiction with the minimality of 
$(|\underline{\lambda}|,|\underline{\mu}|)$. We are not allowed however 
to use braid moves of the types $\tau_{ik}'$ and $\tau_{ik}''$ since 
these change $\lambda_{k}$ while we want to preserve at every braid 
move the equality $\lambda_{1}=\cdots =\lambda_{g}=1$. We thus have the following: $G_1=\langle t_{1},\ldots ,t_{n} \rangle$ 
has orbits $\Sigma_1=\{a_1,\ldots,a_{d-1}\}$, $\Sigma_2=\{a_d\}$, 
$\mu_{k}=(a_ca_d)$, $G_2=\langle \mu_{k},\ldots,\mu_{g}\rangle \subseteq 
S(\{a_c,a_d\})\times S(\Sigma_1-\{a_c\})$. Here if $\ell>k$ we have 
either $\mu_{\ell}=\mu_{\ell}'\in S(\Sigma_1-\{a_c\})$ or 
$\mu_{\ell}=(a_ca_d)\mu_{\ell}''$ with $\mu_{\ell}''\in 
S(\Sigma_1-\{a_c\})$. If $\mu_{\ell}'\neq 1$ or $\mu_{\ell}''\neq 
1$ respectively we might apply Mochizuki's Lemma~\ref{s4.37a} and the braid move 
$\rho_{1\ell}'$ in order to decrease $|\mu_{\ell}|$. This is 
impossible by the minimality of 
$(|\underline{\lambda}|,|\underline{\mu}|)$ so for every $\ell >k$ 
either $\mu_{\ell}=1$ or $\mu_{\ell}=(a_ca_d)$
\end{stepc}
\begin{stepc}\label{c3}
\emph{We claim
$\mu_{\ell}=1$ for $\forall \ell <g$,\; $\mu_{g}=(a_ca_d)$
}.\quad 
Since renumbering and braid moves are commutative operations we may 
assume without loss of generality that $\Sigma_1=\{1,\ldots,d-1\}$ and 
$c=1$ so every $\mu_{\ell}$ equals either 1 or $(1d)$. Let $k$ be the 
minimal index such that $\mu_{k}=(1d)$. If $k=g$ there is nothing to 
prove, so suppose that $k<g$. Replacing eventually $(t_1,\ldots,t_n)$ by a 
braid-equivalent sequence we may assume that 
$(t_1,\ldots,t_n;\lambda_1,\mu_1,\ldots,\lambda_g,\mu_g)$ has the 
following form:
\begin{equation}\label{es5.70bis}
((12),(12),(13),(13),\ldots,(1\, d-1)(1\, 
d-1);1,1,\ldots,\underset{k}{1,(1d)},\ldots).
\end{equation}
First suppose that $\mu_{g}=1$. We claim that replacing $\mu_{g}$ by $(1d)$ 
one obtains a braid-equivalent Hurwitz system. Indeed $\tau_{1k}'$ 
(cf. Corollary~\ref{s2.30a}) transforms Eq.\eqref{es5.70bis} into 
$((2d),(12),\ldots ;\linebreak
\ldots ,(2d),(1d),\ldots )$. Apply the following 
sequence of braid-equivalences: replace $(2d),(12)$ by $(12),(1d)$; 
move this pair to the end by braid moves of type $\sigma_{j}'$; perform 
$\rho_{ng}''$; move the pair $(12),(1d)$ backward to  the front by 
braid moves of type $\sigma_{j}''$; replace it by $(2d),(12)$ and finally 
apply $\tau_{1k}''$. One obtains 
\begin{equation}\label{es5.70abis}
((12),(12),(13),(13),\ldots,(1\, d-1)(1\, 
d-1);1,1,\ldots,\underset{k}{1,(1d)},\ldots,1,(1d))
\end{equation}
as claimed. We may thus suppose that Eq.\eqref{es5.70bis} equals 
Eq.\eqref{es5.70abis}. Applying $\tau_{1g}'$ to Eq.\eqref{es5.70abis} and 
replacing $(2d),(12)$ by $(1d),(2d)$ we obtain
\[
((1d),(2d),(13),(13),\ldots ;\ldots ,1,(1d),\ldots,(2d),(1d))
\]
Applying $\rho_{1k}''$ one transforms $\mu_{k}'=(1d)$ into 1 without 
changing $t_1=(1d)$. Replacing $(1d),(2d)$ by $(2d),(12)$ and applying 
$\tau_{1g}''$ we obtain a Hurwitz system with $\lambda_{1}'=\cdots 
=\lambda_{g}'=1$, $\mu_{1}'=\cdots =\mu_{k}'=1$ and $\mu_{\ell}'=1$ or 
$(1d)$ for $\ell >k$. This shows that the assumption $\mu_{k}\neq 1$ 
for $k<g$ contradicts the minimality of 
$(|\underline{\lambda}|,|\underline{\mu}|)$.
\end{stepc}
We have so far worked with Hurwitz systems. When working with 
equivalence classes  we may moreover renumber arbitrarily $\{1,\ldots 
,d\}$. We thus conclude that the equivalence class (modulo inner 
automorphisms) of a Hurwitz system of the type of Step~\ref{c3} is 
braid-equivalent to the equivalence class Eq.\eqref{es5.70bis} with 
$k=g$, so it has the normal form Eq.\eqref{es5.65bis} as claimed
\end{proof}

\section{The case $n=2d-6,\: g=1$}\label{s6}
Unless otherwise specified in this section we shall work with Hurwitz 
systems $(t_1,\ldots,t_n;\lambda,\mu)$ where $n\geq 2$ and 
$t_{1},\ldots ,t_{n}$ are transpositions in $S_d$.
\begin{pro}\label{s6.72}
Let $(t_1,\ldots,t_n;\lambda,\mu)$, with $\lambda=1$, be a Hurwitz 
system with monodromy group $S_d$. Suppose that $d-1\leq 
n<2d-2$. Let $e=\frac{n}{2}+1$. Then the equivalence class 
$[t_1,\ldots,t_n;\lambda,\mu]$ is braid-equivalent to 
\begin{equation}\label{es6.72}
[(12),(12),(13),(13),\ldots,(1e),(1e);1,(1\, e+1 \ldots d)]
\end{equation}
\end{pro}
\begin{proof}
We may assume without loss of generality that $(|\lambda|,|\mu|)$ is 
minimal and $\ell(\underline{\Sigma})$ is maximal in lexicographic 
order as in \ref{s5.63bbis}~(*). By Lemma~\ref{s5.61bis} the subgroup 
$\langle t_{1},\ldots ,t_{n} \rangle \subseteq S_d$ has a single 
nontrivial domain of transitivity $\Sigma_1$ and furthermore since 
$t_{1}\cdots t_{n}=1$ we may assume that $t_i= t_{i+1}$ for $i\equiv 
1(mod\, 2)$. Let $\Sigma_1=\{a_1,\ldots,a_d\}$. Decomposing $\mu$ into 
a product of independent cycles $\mu = \mu_1\cdots \mu_k$ we have  
each $\mu_i$ contains at most one element of $\Sigma_1$. Indeed, 
otherwise using Mochizuki's Lemma~\ref{s4.37a} and applying the braid 
move $\rho_{ng}''$ (cf. Eq.\eqref{es2.30b}) we might decrease $|\mu|$. 
Furthermore the transitivity of $\langle t_{1},\ldots ,t_{n},\mu 
\rangle$ implies that every $\mu_i$ contains exactly one element of 
$\Sigma_1$ and each element of $\overline{\Sigma}_1$ is contained in 
one of $\mu_{i}$. The inequality $n\geq d-1$, equivalent to $e>d-e$, 
implies that at least one element of $\Sigma_1$ does not appear in the 
 cycles $\mu_{i},\; 1\leq i\leq k$. We claim that if $k\geq 
2$ then there is a braid-equivalent Hurwitz system 
$(t'_1,\ldots,t'_n;1,\mu')$ such that $\mu'$ is a cycle of order 
$d-e+1$ containing a single element of $\Sigma_1$. Braid moves and
renumbering of $\{1,\ldots ,d\}$ are commutative operations, so for 
proving the claim we may assume without loss of generality that 
$(t_1,\ldots,t_n;1,\mu)$ equals 
\begin{equation}\label{es6.74a}
((12),(12),\ldots,(1e),(1e);1,(2\; e+1\ldots e+i_1)(3\; e+i_1+1 \ldots 
)\ldots ).
\end{equation}
Performing $\tau_{11}'$ we obtain 
\begin{equation*}
((1\, e+1),(12),\ldots,(1e),(1e);(1\, e+1),(2\, 
e+1\ldots)(3\ldots)\ldots).
\end{equation*}
Moving the pair $(1\, e+1),(12)$ to the end, replacing it by 
$(12),(2\, e+1)$ and performing $\rho_{n1}''$ we obtain 
\begin{equation}\label{es6.74b}
((13),(13),\ldots,(1e),(1e),(12),(12);(1\, e+1),(2\, e+2\ldots)(3\, 
e+i_1+1\ldots)\ldots).
\end{equation}
Let us consider the system
\begin{equation}\label{es6.74c}
((12),(12),(13),(13),\ldots,(1e),(1e);1,(2\; e+2\ldots e+i_1)(3\; 
e+1\; e+i_1+1\ldots)\ldots).
\end{equation}
Moving the pair $(13),(13)$ to the front and then performing 
$\tau_{11}'$ we obtain
\[
((1\, e+1),(13),(12),(12),(14),(14),\ldots ;(1\, e+1),(2\, 
e+2\ldots)(3\, e+1\ldots)\ldots).
\]
Moving $(1\, e+1),(13)$ to the end, replacing it by $(13),(3\, e+1)$ 
and performing $\rho_{n1}''$ we obtain 
\begin{equation}\label{es6.75a}
((12),(12),(14),\ldots,(1e),(13),(13);(1\; e+1),(2\; 
e+2\ldots)(3\; e+i_1+1\ldots)\ldots).
\end{equation}
Since Eq.\eqref{es6.74c} is braid-equivalent to Eq.\eqref{es6.75a} and 
Eq.\eqref{es6.75a} is braid-equivalent to Eq.\eqref{es6.74b} we 
conclude that Eq.\eqref{es6.74a} is braid-equivalent to Eq.\eqref{es6.74c}. 
The sequence Eq.\eqref{es6.74c} is obtained from Eq.\eqref{es6.74a} by 
moving one element from the cycle $\mu_{1}$ to the cycle $\mu_{2}$.  Repeating 
this transformation several times we  obtain a cycle $\mu'$ of order $d-e+1$. 
We thus proved that 
replacing the given Hurwitz system by a braid-equivalent one we may 
assume that $\langle t_{1},\ldots ,t_{n} \rangle$ has a single nontrivial 
orbit $\Sigma_1=\{a_1,\ldots,a_e\}$ and $\mu=(a_ia_{e+1}\ldots a_d)$ 
where $1\leq i\leq e$ and $\{a_1,a_2,\ldots,a_d\}=\{1,\ldots ,d\}$. 
Renumbering and normalizing the sequence 
$(t_1,\ldots,t_n)$ according to the classical result of Clebsch and 
Hurwitz we obtain that $[t_1,\ldots,t_n;1,\mu]$ is braid-equivalent to 
the equivalence class Eq.\eqref{es6.72}. 
\end{proof}
\begin{rem}\label{s6.76a}
The proof of Proposition~\ref{s6.72} may be easily adapted in order to 
show that  every equivalence class 
$[t_1,\ldots,t_n;\lambda,1]$ with full monodromy group $S_d$ and $
d-1\leq n <2d-2$ is braid-equivalent to 
\begin{equation}\label{es6.76}
[(12),(12),\ldots,(1e),(1e);(1\, e+1\ldots d),1].
\end{equation}
Here one should  modify Condition~\ref{s5.63bbis}(*) assuming 
without loss of generality $(|\lambda|,|\mu|)=(|\lambda|,0)$ is minimal 
in the reverse lexicographical order and $\ell(\underline{\Sigma})$ is 
maximal in lexicographical order. One then repeats the proof of the 
proposition using $\rho_{n1}'$ instead of $\tau_{11}'$ and 
$\tau_{11}''$ instead of $\rho_{n1}''$.
\end{rem}
\begin{lem}\label{s6.76}
Let $(t_1,\ldots,t_n;\lambda,\mu)$ be a Hurwitz system with monodromy 
group $S_d$. Let $\underline{\Sigma}=\{\Sigma_1,\ldots,\Sigma_m\}$  be 
the domains of transitivity of $\langle t_{1},\ldots ,t_{n} \rangle$ 
ordered in such a way that $\# \, \Sigma_1\geq \# \, \Sigma_2\geq 
\cdots \geq\# \, \Sigma_m$. Suppose that $(|\lambda|,|\mu|)$ is minimal and 
$\ell(\underline{\Sigma})$ is maximal as in \ref{s5.63bbis}(*). 
Suppose that $\lambda\neq 1$. 
Then the following three conditions cannot hold simultaneously.
\renewcommand{\theenumi}{\roman{enumi}}
\begin{enumerate}
\item
$\Sigma_1$ is non-minimal (cf. \ref{s5.60bis}).
\item There exists $a\in \Sigma_1$ such that $a^{\lambda}\notin 
\Sigma_1$ (equivalently there exists $c\in \Sigma_1$ such that 
$c^{\lambda^{-1}}\notin \Sigma_1$).
\item There exists $b\in \Sigma_1$ such that $b^{\lambda}=b$.
\end{enumerate}
\end{lem}
\begin{proof}
Decomposing $\lambda$ into a product of independent cycles 
$\lambda=\lambda_1\cdots \lambda_k$ we see that each of the two 
conditions of (ii) is equivalent to the existence  of a $\lambda_i$ 
which contains elements of both $\Sigma_1$ and 
$\overline{\Sigma}_1=\{1,\ldots ,d\}-\Sigma_1$. Let $\mu=\mu_1\cdots 
\mu_{\ell}$ be the decomposition of $\mu$ into a product of 
independent cycles. The minimality of $(|\lambda|,|\mu|)$ implies that neither
 $\lambda_i$ nor $\mu_j$ may contain two different elements of some 
$\Sigma_r$. In fact if this were the case then we might apply 
Mochizuki's Lemma~\ref{s4.37a} and one of the braid moves 
$\tau_{11}''$ or $\rho_{n1}''$ in order to decrease 
$(|\lambda|,|\mu|)$. In particular for each $x\in \Sigma_1$ we have an 
alternative: either 
$x^{\lambda}=x$ or $x^{\lambda}\notin \Sigma_1$.

Suppose, by way of contradiction, that (i), (ii) and (iii) hold simultaneously.
 Let $s=[\lambda,\mu]$. 
By Lemma~\ref{s5.61bis} we must have $n+|s|<2d-2$ since otherwise 
$\langle t_{1},\ldots ,t_{n} \rangle$ is transitive, so applying 
Mochizuki's Lemma~\ref{s4.37a} and the braid move $\tau_{11}''$ we 
might decrease $|\lambda|$ which contradicts the minimality of 
$(|\lambda|,|\mu|)$ assumed in the lemma. Lemma~\ref{s5.61bis} yields 
moreover that $n_1+|s|_{\Sigma_1}|=2(\# \, \Sigma_1 -1)$, so 
Condition~(i) means $s|_{\Sigma_1}$ has $q\geq 2$
domains of transitivity $\Delta_{1},\ldots,\Delta_{q}$. We claim one 
can choose $a$ and $b$ satisfying (ii) and (iii) respectively so that 
they belong to different domains. In fact if $b\in \Delta_{k}$ and all 
$x\in \Sigma_1$ with $x^{\lambda}\notin \Sigma_1$ belong to 
$\Delta_{k}$ then we may replace $b$ by an arbitrary element of some 
other $\Delta_{\ell}$ according to the alternative of the preceding 
paragraph. Let $a\in \Delta_{i},\; b\in \Delta_{j},\; 
i\neq j$. Let $T_1$ be the subsequence of $T=(t_1,\ldots,t_n)$ 
consisting of the transpositions that move elements of $\Sigma_1$. Let 
$a'=a^{s},\; b'=b^{s}$. Using Corollary~\ref{s4.41} we can replace 
$T_1$ by a sequence $(\cdots,(a'b'),(a'b'))$ and then move the pair 
$(a'b'),(a'b')$ to the end of $T$. We then perform the braid move 
$\rho_{n1}'$ which transforms $t_n=(a'b')$ into 
$(a'b')^{s^{-1}\lambda} = (ab)^{\lambda} = (a^{\lambda}b)$. The group 
generated by the new sequence $\langle T'\rangle$ has orbit 
$\Sigma_1\cup \{a^{\lambda}\}$. By Mochizuki's Lemma~\ref{s4.37a} one 
may replace $T'$ by $((a\, a^{\lambda}),\ldots)$. Performing 
$\tau_{11}''$ one transforms $\lambda$ into $(a\, a^{\lambda})\lambda$ 
and clearly $|(a\, a^{\lambda})\lambda|<|\lambda|$. This contradicts 
the minimality of $(|\lambda|,|\mu|)$.
\end{proof}
\begin{lem}\label{s6.79}
Let $(t_1,\ldots,t_n;\lambda,\mu)$ be a Hurwitz system with monodromy 
group $S_d$. Let $d-1\leq n<2d-2$. Suppose that $(|\lambda|,|\mu|)$ is 
minimal and $\ell(\underline{\Sigma})$ is maximal as in 
\ref{s5.63bbis}(*). Suppose that $s=[\lambda,\mu]=1$. Then $\lambda=1$ and the equivalence class $[t_1,\ldots,t_n;\lambda,\mu]$ is 
braid-equivalent to Eq.\eqref{es6.72}.
\end{lem}
\begin{proof}
Let $\underline{\Sigma}=\{\Sigma_1,\ldots,\Sigma_m\}$ be the orbits of 
$\langle t_{1},\ldots ,t_{n} \rangle$. Let $e=\frac{n}{2}+1$. 
According to Lemma~\ref{s5.61bis} we have $\# \, \Sigma_i=1$ for each 
$i\geq 2$ and $\# \, \Sigma_1=e$. Suppose, by way of contradiction, that 
$\lambda\neq 1$. According to Lemma~\ref{s6.76} its conditions (i) -- (iii)  
cannot hold simultaneously. The orbit $\Sigma_1$ is non-minimal since 
$s|_{\Sigma_1}=1$. As we saw in the course of the proof of 
Lemma~\ref{s6.76} the minimality of 
$(|\lambda|,|\mu|)$ implies that if one decomposes $\lambda$ into a 
product of independent cycles $\lambda=\prod\lambda_i$ then every 
$\lambda_i$ contains at most one element of $\Sigma_1$. The inequality 
$n\geq d-1$ implies that $\# \, \Sigma_1=e>d-e$. Thus there exists an 
element $b\in \Sigma_1$ such that $b^{\lambda}=b$. Therefore 
Condition~(ii) of Lemma~\ref{s6.76} must fail to hold, i.e. we have
$a^{\lambda}=a$ for every $a\in \Sigma_1$. Let $F\subseteq \{1,\ldots 
,d\}$ be the set of fixed points of $\lambda$. We have $F\supset 
\Sigma_1$ so the group $\langle t_{1},\ldots ,t_{n} \rangle$ leaves 
$F$ invariant. Since $[\lambda,\mu]=1$ we have  $F^{\mu}=F$. 
Therefore the assumption $\lambda\neq 1$ implies that $\langle t_{1},\ldots 
,t_{n},\lambda,\mu \rangle$ is not a transitive group which is a 
contradiction. The last part of the lemma refers to 
Proposition~\ref{s6.72}.
\end{proof}
\begin{rem}\label{s6.80}
This lemma may be stated differently: if $[\lambda,\mu]=1$ and 
$\lambda\neq 1$ then one can decrease $|\lambda|$ by a sequence of 
braid moves (at the expense of possible increasing of $|\mu|$). 
Modifying appropriately Lemma~\ref{s6.76} one can prove similarly that 
given a Hurwitz system $(t_1,\ldots,t_n;\lambda,\mu)$ with full 
monodromy group $S_d$ and such that $d-1\leq n<2d-2$, then 
$[\lambda,\mu]=1$ and $\mu\neq 1$ imply that one can decrease $|\mu|$ by a 
sequence of braid moves (at the expense of possible increasing of 
$|\lambda|$).
\end{rem}
\begin{thm}\label{s6.81}
Let $d\geq 4,\; n=2d-6$. Let $Y$ be an elliptic curve. Then the 
Hurwitz space $\mathcal{H}^0_{d,n}(Y)$ parametrizing simple coverings 
with monodromy groups equal to $S_d$ is irreducible.
\end{thm}
\begin{proof}
By Theorem~\ref{s2.30c} it suffices to show that given a Hurwitz 
system $(t_1,\ldots,t_n;\lambda,\mu)$ with monodromy group $S_d$ its 
equivalence class is braid-equivalent to
\begin{equation}\label{es6.81}
[(12),(12),(13),(13),\ldots,(1\, d-2),(1\, d-2);1,(1\, d-1\, d)].
\end{equation}
Without loss of generality we may assume that 
$(t_1,\ldots,t_n;\lambda,\mu)$ has minimal $(|\lambda|,|\mu|)$ and 
maximal $\ell(\underline{\Sigma})$ as in \ref{s5.63bbis}(*). Here 
$\underline{\Sigma}=\{\Sigma_1,\ldots,\Sigma_m\}$ are the domains of 
transitivity of $G_1=\langle t_{1},\ldots ,t_{n} \rangle$ ordered in 
such a way that 
$\# \, \Sigma_1\geq \# \, \Sigma_2\geq 
\cdots \geq \# \, \Sigma_m$. We aim to prove that $\lambda=1$ in order 
to apply Proposition~\ref{s6.72}. Let $s=[\lambda,\mu]$. Our first 
goal is to prove $s=1$. Suppose, by way of contradiction, that $s\neq 1$, so 
in particular $\lambda\neq 1$ and $\mu\neq 1$. Since $[\lambda,\mu]$ 
is an even permutation we have $|s|\equiv 0(mod\: 2)$. If $s\geq 4$ 
then $n+|s|\geq 2d-2$, so by Lemma~\ref{s5.61bis} the group 
$G_1=\langle t_{1},\ldots ,t_{n} \rangle$ is transitive. Using 
Mochizuki's Lemma~\ref{s4.37a} and the braid move $\tau_{11}''$ we can 
transform $\lambda$ into $\lambda'$ with $|\lambda'|<|\lambda|$. This 
contradicts the minimality of $(|\lambda|,|\mu|)$. Thus if $s\neq 1$ 
then $|s|=2$. By Lemma~\ref{s5.61bis} the group $G_1$ has two domains 
of transitivity $\underline{\Sigma}=\{\Sigma_1,\Sigma_2\}$ and 
$\Sigma_2$ is minimal. If $\lambda=\lambda_1\cdots 
\lambda_k,\; \mu=\mu_1\cdots \mu_{\ell}$ are the factorizations into 
independent cycles, then as we saw in the proof of Lemma~\ref{s6.76} 
the minimality of $(|\lambda|,|\mu|)$ implies that each $\lambda_i$ or 
$\mu_j$ is a transposition of the type $(ac)$ where $a\in \Sigma_1,\; 
c\in \Sigma_2$. It is clear that $\max\{k,\ell\}\leq \# \, \Sigma_2$.
Let $s_1=(s|_{\Sigma_1})^{\sim},\; s_2=(s|_{\Sigma_2})^{\sim}$ where 
${}^{\sim}$ 
means trivial extension from $S(\Sigma_i)$ to $S_d$. We have 
$s=s_1s_2$ and the following cases may occur: a) $|s_1|=2,\; |s_2|=0$;\; 
b)\; $|s_1|=1,\; |s_2|=1$;\; c)\; $|s_1|=0,\; |s_2|=2$.

In Case (a) by minimality $\# \, \Sigma_2=1$. So $\lambda_1=(a_1c),\; 
\mu=(a_2c)$ where $\{c\}=\Sigma_2$ and $a_1\neq a_2$ since 
$[\lambda,\mu]\neq 1$. Then $[\lambda,\mu]=(a_2a_1c)$. This 
contradicts $[\lambda,\mu]=t_{1}\cdots t_{n}\in S(\Sigma_1)\times 
S(\Sigma_2)$, so Case (a) is impossible.

In Case (b) we have $\# \, \Sigma_2=2$ by minimality and $\# \, 
\Sigma_1\geq 2$. If $d\geq 5$ then $\# \, \Sigma_1\geq 3$. So 
$\Sigma_1$ is non-minimal and furthermore there exists an element $b\in 
\Sigma_1$ such that $b^{\lambda}=b$. If $\lambda=\lambda_1$ or 
$\lambda=\lambda_1\lambda_2$ with $\lambda_1=(ac)$ then $a\in 
\Sigma_1,\; a^{\lambda}\notin \Sigma_1$. This contradicts 
Lemma~\ref{s6.76}. It remains to consider the case $d=4$. Let 
$\Sigma_1=\{a_1,a_2\},\; \Sigma_2=\{c_1,c_2\}$. 

Suppose that $\lambda$ is a 
transposition. We may assume that $\lambda=(a_1c_1)$. We exclude case by 
case the possibilities for $\mu$. The cases $\mu=(a_1c_1)$, or 
$(a_2c_2)$, or $\mu=(a_1c_1)(a_2c_2)$ are impossible since 
$[\lambda,\mu]\neq 1$. The cases $\mu=(a_1c_2)$ or $\mu=(a_2c_1)$ are 
impossible either since $[\lambda,\mu]=s_1s_2$ is not a cycle of order 
3. It remains to exclude the possibility $\mu=(a_1c_2)(a_2c_1)$. Indeed here
\[
[\lambda,\mu] = (a_1c_1)(a_1c_2)(a_2c_1)(a_1c_1)(a_1c_2)(a_2c_1)= 
(a_1c_1)(a_2c_2),
\]
while $[\lambda,\mu]$ must equal $(a_1a_2)(c_1c_2)$. 

Suppose that $\lambda=\lambda_1\lambda_2$. We may assume that $\lambda = 
(a_1c_1)(a_2c_2)$. Then $[\lambda,\mu]$ might be $\neq 1$ only if 
$\mu=(a_1c_2)$ or $\mu=(a_2c_1)$. The calculation above shows that in these 
cases $[\lambda,\mu]\neq (a_1a_2)(c_1c_2)$
So Case (b) is impossible.

In Case (c) we have $\# \, \Sigma_2=3$ by minimality and 
$s|_{\Sigma_1}=1$. Here $d=\# \, \Sigma_1 + \# \, \Sigma_2 \geq 6$. If 
$d\geq 7$ we obtain a contradiction with Lemma~\ref{s6.76} by the same 
argument we used in Case (b). The cases $d=6$ and $\lambda=\lambda_1$ 
or $\lambda=\lambda_1\lambda_2$ are excluded in the same way. It 
remains to consider the case $\lambda=(a_1c_1)(a_2c_2)(a_3c_3)$ where 
$\Sigma_1=\{a_1,a_2,a_3\},\; \Sigma_2=\{c_1,c_2,c_3\}$. Let $T_i$ be 
the subsequence of $T=(t_1,\ldots,t_6)$ which moves elements of 
$\Sigma_i,\; i=1,2$. Replacing $T$ by a braid-equivalent sequence we 
may assume that $T=T_1T_2$ where $T_1=(a_1a_2)(a_1a_2)(a_1a_3)(a_1a_3)$. By 
Main Lemma~\ref{s3.32} the Hurwitz system $(T_1T_2;\lambda,\mu)$ is 
braid-equivalent to $({T_1}^{\mu}T_2;\lambda,\mu)$. If 
$\#({\Sigma_1}^{\mu}\cap \Sigma_2)<3$ then the group $\langle 
{T_1}^{\mu},T_2\rangle$ has a domain of transitivity 
$\Sigma'_1={\Sigma_1}^{\mu}\cup \Sigma_2$ and 
$\#\Sigma'_1>\#\Sigma_1$. This contradicts the maximality of 
$\ell(\underline{\Sigma})$. Hence ${\Sigma_1}^{\mu}=\Sigma_2,\: 
{\Sigma_2}^{\mu}=\Sigma_1$. This implies that the group 
$(t_1,\ldots,t_n;\lambda,\mu)$ is imprimitive which is excluded by 
hypothesis. We thus excluded all possible cases with $s\neq 1$. 
Therefore $s=[\lambda,\mu]=1$.

Suppose that $d\geq 5$. Then from Lemma~\ref{s6.79} and 
Proposition~\ref{s6.72} it follows that $\lambda=1$ and 
$[t_1,\ldots,t_n;\lambda,\mu]$ is braid-equivalent to Eq.\eqref{es6.81} 
which proves the theorem when $d\geq 5$. So, the only case that 
remains to be considered is $d=4, n=2, [\lambda,\mu]=1$. Here we have 
$\underline{\Sigma}=\{\Sigma_1,\Sigma_2,\Sigma_3\}$ with 
$\#\Sigma_1=2,\; \# \, \Sigma_2=\# \, \Sigma_3=1$. Let 
$\overline{\Sigma}_1=\Sigma_2\cup \Sigma_3$. Suppose, by way of 
contradiction, that $\lambda\neq 1$. Unless $\lambda=(ac)(bd)$ where 
$\Sigma_1=\{a,b\},\; \overline{\Sigma}_1=\{c,d\}$ we may apply the 
argument of Lemma~\ref{s6.79} and obtain a contradiction with the 
minimality of $(|\lambda|,|\mu|)$. We obtain a Hurwitz system of the 
type $((ab),(ab);(ac)(bd),\mu)$ where $\mu$ commutes with $(ac)(bd)$. 
The centralizer of $\lambda=(ac)(bd)$ is 
$\{1,(ac)(bd),(ad)(bc),(ab)(cd),(ac),(bd),(abcd),(adcb)\}$. The braid 
move $\rho_{21}''$ transforms $\mu$ into $\mu''=(ab)\mu$ (cf. 
Eq.\eqref{es2.30b}). Thus the minimality of $(|\lambda|,|\mu|)$ excludes 
the cases $\mu=(ab)(cd), (abcd)$ or $(adcb)$. The hypothesis of 
primitivity of $\langle t_1,t_2,\lambda,\mu\rangle$ excludes the cases 
$\mu=1,(ac)(bd)$ or $(ad)(bc)$. It thus remains to consider the cases 
$\mu=(ac)$ or $\mu=(bd)$. These are equivalent up to reordering, so it 
suffices to consider $\mu=(ac)$. We perform the following sequence of 
braid moves starting from $((ab),(ab);(ac)(bd),(ac))$. Applying 
$\tau_{11}'$ we obtain $((bc),(ab);(acdb),(ac))$. Replacing 
$(bc),(ab)$ by $(ac),(bc)$ and then applying $\tau_{11}''$ (cf. 
Theorem~\ref{s2.27}~(f)) we obtain $(ac)\mapsto (ac)^{(abc)}=(ab)$ 
and $(acdb)\mapsto (ac)(acdb)=(adb)$. So we obtain 
$((ab),(bc);(bad),(ac))$. Applying $\tau_{11}''$ again we transform 
$(\ldots; (bad),(ac))$ into $(\ldots;(bd),(ac))$. This contradicts the 
minimality of $(|\lambda|,|\mu|)$ of the initial Hurwitz system 
$((ab),(ab);(ac)(bd),(ac))$. All possible cases with 
$[\lambda,\mu]=1,\; \lambda\neq 1$ were excluded, so $\lambda=1$. We 
conclude that $(t_1,t_2;\lambda,\mu)=((ab),(ab);1,\mu)$. The minimality of 
$(|\lambda|,|\mu|)=(0,|\mu|)$ implies that no independent cycle of the 
factorization of $\mu$ may contain both $a$ and $b$. Since $\langle 
(ab),\mu\rangle$ is transitive we conclude that the factorization is either 
$\mu=\mu_{1}\mu_{2}=(ac)(bd)$ with $\{b,d\}=\overline{\Sigma}_1$ or 
$\mu$ is a cycle of order 3 containing the two elements of 
$\overline{\Sigma}_1$. The former case is impossible since $\langle 
(ab),\mu\rangle$ is a primitive group. We conclude that the equivalence 
class $[(ab),(ab);1,\mu]$ equals $[(12),(12);1,(134)]$. The theorem is proved.
\end{proof}
We can sharpen Theorem~\ref{s4.55} using the same proof and 
Theorem~\ref{s4.45},\; Theorem~\ref{s5.64bis} and Theorem~\ref{s6.81}.
\begin{thm}\label{s6.87}
Let $Y$ be a smooth, projective curve of genus $\geq 1$. 
Let $d\geq 2,\; n>0,\; n\equiv 0(mod\; 2)$. If $d$ is prime then 
$\mathcal{H}_{d,n}(Y) = \mathcal{H}^0_{d,n}(Y)$. If $d$ is not prime
let $d'$ be its maximal divisor  $\neq d$. 
 Suppose that $n\geq \max(2,2d'-4)$ or if 
$g=1$ suppose that $n\geq \max(2,2d'-6)$.
 Then the connected 
components of $\mathcal{H}_{d,n}(Y)$ which correspond to simple coverings 
with monodromy groups $\neq S_d$ are in one-to-one correspondence with 
the equivalence classes of unramified coverings $[\tilde{Y}\to Y]$ 
of degrees $d_2|d$ where $d_2\neq 1, d$.
\end{thm}

\bibliographystyle{amsalpha}
\providecommand{\bysame}{\leavevmode\hbox to3em{\hrulefill}\thinspace}

\bigskip

\noindent
{\sc
Dipartimento di Matematica, Universit\`{a} di Palermo\\
Via Archirafi n.34, 90123 Palermo, Italy}\qquad and 
\par \smallskip \noindent
{\sc Institute of Mathematics, Bulgarian Academy of Sciences}

\medskip \noindent
{\it E-mail address:} {\bf kanev@math.unipa.it}

\end{document}